\documentclass[a4paper,twoside,reqno,11pt]{amsart}

\usepackage{a4wide,verbatim}
\usepackage[english]{babel}
\usepackage{epsfig,rotating,color}
\usepackage{amsmath,amsfonts,amssymb,amsgen,amsbsy}
\usepackage{eucal,esint,dsfont,mathrsfs,MnSymbol}
\usepackage{stmaryrd} 
\usepackage{multicol}
\usepackage{tikz,pgfplots}
\usepackage{graphicx}
\usepackage{mathtools}

\numberwithin{equation}{section}

\usepackage{caption}

\usepackage[square,numbers]{natbib}
\usepackage{datetime} 
\newif\ifdraft

\ifdraft
\usepackage[notref,notcite]{showkeys}
\fi

\usepackage{hyperref}

\hypersetup{
colorlinks=true, 
breaklinks=true, 
urlcolor= blue, 
linkcolor= blue, 
bookmarksopen=true, 
pdftitle={}, 
pdfauthor={}, 
pdfsubject={} 
}

\newcommand{\R}{{\mathds R}}

\newcommand{\eqdef}{\stackrel{\text{\tiny{def}}}{=}} 

\newcommand{\ud}{\mathrm{d}}

\newcommand{\half}{{\textstyle{1\over2}}}

\newcommand{\third}{{\textstyle{1\over3}}}
\newcommand{\quat}{{\textstyle{1\over4}}}

\newcommand{\sixth}{{\textstyle{1\over6}}}
\newcommand{\twothird}{{\textstyle{2\over3}}}
\newcommand{\threehalf}{{\textstyle{3\over2}}}

\newcommand{\sgn}{\operatorname{sgn}}

\newlength{\intwidth}

\ifdraft
  \title[DRAFT \shortdate\today\quad \currenttime]
  {\bf Global weak solutions of a Hamiltonian regularised Burgers equation}
  \author[DRAFT \shortdate\today\quad \currenttime]
  {Billel GUELMAME, St\'ephane JUNCA, Didier CLAMOND and Robert L. PEGO}
\else
  \title[Regularised Burgers equation]{\bf Global weak solutions of a Hamiltonian regularised Burgers equation}
  \author[Guelmame et al.]{Billel GUELMAME, St\'ephane JUNCA, Didier CLAMOND and Robert L. PEGO}
\fi

\newcommand{\nfont}{\fontshape{n}\selectfont}
\address{({\nfont\textbf{Billel Guelmame}}) Universit\'e C\^ote d'Azur, CNRS, Inria,  LJAD, France.} 
\email{billel.guelmame@univ-cotedazur.fr}
\address{({\nfont\textbf{St\'ephane Junca}}) Universit\'e C\^ote d'Azur, CNRS,  Inria, LJAD, France.} 
\email{stephane.junca@univ-cotedazur.fr}
\address{({\nfont\textbf{Didier Clamond}}) Universit\'e C\^ote d'Azur, CNRS, LJAD, France.} 
\email{didier.clamond@univ-cotedazur.fr}
\address{({\nfont\textbf{Robert L. Pego}}) Department of Mathematical Sciences and Center for Nonlinear Analysis, Carnegie
Mellon University, Pittsburgh, Pennsylvania, PA 12513, USA.} 
\email{rpego@cmu.edu}

\date{\today}

\usepackage{amsthm}
\newtheorem{thm}{Theorem}[section]
\newtheorem{definition}[thm]{Definition}
\newtheorem{lem}[thm]{Lemma}
\newtheorem{pro}[thm]{Proposition}

\theoremstyle{remark}
\newtheorem{rem}[thm]{Remark}

\usepackage{enumitem}
\newlist{steps}{enumerate}{1}
\setlist[steps, 1]{label = Step \arabic*:}

\graphicspath{ {images/} }

\usepackage{empheq,etoolbox}
\patchcmd{\subequations}
  {\theparentequation\alph{equation}}
  {\theparentequation.\alph{equation}}
  {}{}

\usepgfplotslibrary{fillbetween}
\usetikzlibrary{arrows.meta,intersections}

\pgfmathdeclarefunction{bicuadratic}{1}{%
    \pgfmathparse{tan(40*abs(#1-2.8))+2}%
}

\pgfplotsset{compat=1.15}
\makeatletter
\pgfplotsset{
    mark max/.style={
        point meta rel=per plot,
        visualization depends on={x \as \xvalue},
        scatter/@pre marker code/.code={%
        \ifx\pgfplotspointmeta\pgfplots@metamax
            \def\markopts{mark=none}%
            \coordinate (maximum);
        \fi
            \def\markopts{mark=none}
            \expandafter\scope\expandafter[\markopts]
        },%
        scatter/@post marker code/.code={%
            \endscope
        },
        scatter
    }
}
\newcommand\DrawEpigraph[6][draw=white,top color=gray!80!black!05,bottom color=gray!90!black!80]{
  \coordinate (plot-left) at ([yshift=#4]axis cs:#2,\pgfplots@metamax);
  \coordinate (plot-right) at ([yshift=#5]axis cs:#3,\pgfplots@metamax);
  \path[name path=diagonal,draw=none] (plot-left) -- (plot-right);
  \addplot[#1] fill between[of=#6 and diagonal];
}

\newcommand\Epigraph[6][draw=white,top color=gray!0!black!0,bottom color=gray!0!black!0]{
  \coordinate (plot-left) at ([yshift=#4]axis cs:#2,\pgfplots@metamax);
  \coordinate (plot-right) at ([yshift=#5]axis cs:#3,\pgfplots@metamax);
  \path[name path=diagonal,draw=none] (plot-left) -- (plot-right);
  \addplot[#1] fill between[of=#6 and diagonal];
}

\begin{document}

\begin{abstract} 
A nondispersive, conservative regularisation of the inviscid Burgers equation is proposed and studied.
Inspired by a related regularisation of the shallow water system recently introduced by Clamond and Dutykh,
the new regularisation provides a family of Galilean-invariant interpolants between the inviscid Burgers equation
and the Hunter--Saxton equation. It admits weakly singular regularised shocks and cusped traveling-wave weak solutions.
The breakdown of local smooth solutions is demonstrated, and  
the existence of two types of global weak solutions, conserving or dissipating an $H^1$ energy, is established.
Dissipative solutions satisfy an Oleinik inequality like entropy solutions of the inviscid Burgers equation.
As the regularisation scale parameter $\ell$ tends to $0$ or $\infty$, 
limits of dissipative solutions are shown to satisfy the inviscid Burgers or Hunter--Saxton equation respectively, 
forced by an unknown remaining term.
\end{abstract}


\maketitle
\medskip


\medskip


\tableofcontents

\section{Introduction}\label{Introduction} %


The dispersionless shallow water equations, also called the Saint-Venant equations,
admit shock-wave solutions. Recently, a Hamiltonian regularisation of this system 
(rSV), has been proposed which approximates these discontinuous waves by less
singular ones  \cite{ClamondDutykh2018a}. 
The rSV system can be written 
\begin{subequations}\label{rSV}
\begin{gather}
h_t\ +\,\left[\,h\,u\,\right]_x\ =\ 0, \label{rSV1}  \\
\left[\,h\,u\,\right]_t \ +\, \left[\,h\,u^2\,+\,\half\, g\, h^2\, +\, \varepsilon\, \mathcal{R}\, 
h^2\, \right]_x \ =\ 0, \label{rSV2} \\
\mathcal{R}\ \eqdef\ h\left(\,u_x^{\,2}\,-\,u_{xt}\, -\, u\, u_{xx}\, \right)\, -\ 
g\left(\,h\,h_{xx}\, +\, \half\, h_x^{\,2}\,\right), \label{rSV3} 
\end{gather}
\end{subequations}
where $\varepsilon$ is a small positive parameter, $h$ is the total water depth and $u$ 
is the velocity. The classical Saint-Venant equations can be obtained letting $\varepsilon\to0$.
This regularisation is Galilean invariant, non-dispersive, non-diffusive, 
and conserves energy for regular solutions.  It also admits 
regularised shock-wave weak solutions which have the same wave speed 
and which dissipate energy at the same rate 
as shocks in the classical Saint-Venant (cSV) equations, 
\citep{PuEtAl2018}. 
Some mathematical results on rSV were obtained by \citet{PuEtAl2018} and \citet{LiuEtAl2019}, 
but several natural questions remain open, such as the existence of global weak solutions. 
{Inspired by the rSV equations, a more general regularisation of the unidimensional barotropic Euler system has been derived and studied in \cite{rE2021}}.

In the present work we consider such questions for an analogous but simpler model equation, 
namely a Hamiltonian regularisation of the inviscid Burgers equation $u_t+uu_x=0$.  
Motivated by the rSV and the dispersionless Camassa-Holm \citep{CamassaHolm1993} equations, in 
Section \ref{secheuder} we describe a regularised Burgers equation (rB) in the form
\begin{equation}\label{rB0}
u_t\ +\ u\, u_x\ =\ \ell^2\left( u_{txx}\ +\ 2\, u_x\, u_{xx}\ + u\, u_{xxx} \right),
\end{equation}
where $\ell\geqslant0$ is a parameter. 
Being a scalar equation, the rB equation is more tractable than the rSV system.
An equation mathematically equivalent to \eqref{rB0} has previously appeared in \cite[Remark~1]{DDM} together with a Hamiltonian formulation.
It can be compared to the well-known Camassa--Holm (CH) equation \cite{CamassaHolm1993}, the Degasperis--Procesi (DP) equation \cite{DegasperisProcesi1999} and the Benjamin--Bona--Mahony (BBM) equation  
\cite{BenjaminEtAl1972} (see Section \ref{secheuder} below).

The purpose of the present paper is to establish several basic results for \eqref{rB0}, 
including the existence of local smooth solutions, blow-up, global weak solutions, and
weakly singular traveling waves, and also to study the limiting cases $\ell\to0$ and $\ell\to+\infty$. 

The local (in time) existence of smooth solutions of the rB equation \eqref{rB0} has been established in \cite{Yin2004b,yin2007cauchy} for a generalised Camassa--Holm equation that covers \eqref{rB0} as special case.
The existence of global weak solutions of the Camassa--Holm equation in the space $H^1$ has been widely studied before, we refer to \citep{BressanConstantin2007a,BressanConstantin2007b,ChenTian2009,
CHK,CK2,HR2007,HR2008,HR2009,XinZhang2000}.
We also refer to \cite{CK3,CK4,CK5} for the existence of solutions of the DP equation.


Our treatment of global weak solutions is analogous to the treatment of the Camassa--Holm equation 
by Bressan and Constantin in \cite{BressanConstantin2007a,BressanConstantin2007b}. 
We rewrite (\ref{rB0}) into an equivalent semi-linear system, 
but without asking the initial data to be in $H^1$. 
We then prove the existence of a so-called {\em conservative\/} global weak solution 
(Theorem \ref{existence}, cf.~\cite{BressanConstantin2007a}), which locally conserves energy. 
Energy conservation may not be appropriate for approximating shock waves, however.
We obtain another type of solution called {\em dissipative} 
(see Theorem \ref{existence2} below, cf.~\cite{BressanConstantin2007b}), 
by slightly modifying the equivalent system.
Dissipative solutions satisfy an Oleinik inequality of the form 
\begin{equation}\label{oleinik}
u_x(t,x)\ \leqslant\ \frac Ct\,, \qquad t>0,\ x\in\R.
\end{equation}
This inequality is well known to ensure uniqueness for entropy solutions of the inviscid Burgers equation. 
However, uniqueness for dissipative solutions of rB remains an open problem. 

In order to study the limiting cases $\ell \to 0$ and $\ell \to + \infty$, the equivalent system and the Oleinik inequality \eqref{oleinik} are  used to obtain a uniform  BV 
estimate independent of the parameter $\ell$  for the dissipative solutions (Lemma 
\ref{TV}). When $\ell \to 0$, a dissipative solution converges (up to a subsequence) 
to a function $u$ that satisfies the Burgers equation with a remaining term 
(see Theorem \ref{ThmCompactness} below). 
If the remaining term is zero, then the entropy solution of Burgers is recovered. 
We prove that this term is zero for smooth solutions of Burgers equations (see Proposition \ref{mu}). 
However, the disappearance of the remaining term in general remains an open problem.
Similar results are obtained when $\ell \to +\infty$, where the limit is a solution 
of the Hunter--Saxton equation, at least before the appearance of singularities\footnote{ 
``Singularity'' is used here to describe the blow-up of derivatives, which corresponds 
to shocks of the classical Burgers equation. Contrary to the Burgers case, solutions 
of rB remain continuous at the singularities.}
(Theorem \ref{ThmCompactness2} and Proposition \ref{nu} below).
The limiting case $\ell \to 0$ of the Camassa--Holm equation is more challenging. Indeed, dissipative solutions of the CH equation satisfy an Oleinik inequality with a constant that depends on $\ell$. Thus,
the compactness arguments presented in this paper cannot be used for the CH equation.
However, the limiting case of the viscous CH equation have been studied in \cite{CD2016,CK1,Hwang2007} under the condition ``$\ell$ is small enough compared to the viscosity parameter''. The authors proved that as the viscosity parameter goes to zero, we recover the unique entropy solution of the scalar conservation law $u_t + (3u^2/2)_x =0$.

We find below that the rB equation has a great variety of weakly singular
traveling wave solutions, solutions which are bounded, continuous and piecewise
smooth but which may dissipate (or gain) energy at isolated points where
derivatives become infinite.  
All these waves have analogs for the rSV system \eqref{rSV1}--\eqref{rSV3}.
In particular, corresponding to each simple
shock-wave entropy solution of the inviscid Burgers equation, there is a
monotonic traveling-wave {\em dissipative} solution of the rB equation, 
having the same limiting states, shock speed and energy dissipation rate.
We also find {cusped} traveling waves (both periodic and solitary in nature)
that are {\em conservative}. 
Furthermore, there is a great abundance of composite waveforms that
are neither dissipative nor conservative, 
which were overlooked in \cite{PuEtAl2018} but are similar to some of the
many types of weak traveling wave solutions of the Camassa-Holm equation
found by Lenells~\cite{Lenells2005}.  

This paper is organised as follows. A heuristic derivation of the rB equation is given in section 
\ref{secheuder}. Section \ref{Existence of the solution} is devoted to study the existence of local 
smooth solutions. In Section \ref{Global weak solutions}, a proof of the existence of global 
conservative solutions is given. 
The global dissipative solutions are obtained in Section  \ref{ssecglobexidis}.
Weakly singular traveling waves are described in Section~\ref{s:waves}, 
including energy-conserving `cuspons' and energy-dissipating weakly singular shocks.
Section \ref{Compactness} studies the limiting 
cases $\ell\to0$ and $\ell\to+\infty$ for dissipative solutions. 
The optimality of the requirement that $u_x \in L^2_{loc}$ for weak solutions 
is shown in section \ref{secoptimality}, where we prove in particular 
that when a smooth solution breaks down, 
$u_x$ may blow up in $L^p_{\rm loc}$ for all $p>2$.


\section{Heuristic derivation of a regularised Burgers equation}\label{secheuder}

In order to describe a suitable regularisation 
of the inviscid Burgers equation
with similar features as the rSV system \eqref{rSV},
we note first that the rSV equations yield
\begin{equation}
u_t\ +\ u\,u_x\ +\ g\,h_x\ +\ \varepsilon\left(\,h\,\mathcal{R}_x\,+\,2\,\mathcal{R}\,h_x\,\right)\,=\ 0.
\end{equation}
When $h$ is taken constant, this equation (with the definition of $\mathcal{R}$ given at \eqref{rSV3}) becomes
\begin{equation}\label{h=1}
u_t\ +\ u\,u_x\ =\ \ell^2\left[\,u_{xxt}\, -\,u_x\,u_{xx}\, +\,u\,u_{xxx}\,\right],
\end{equation}
where $\ell\eqdef h\sqrt{\varepsilon}\geqslant0$ is a constant characterising a length scale 
for the regularisation. 
After the change of independent variables $(t,x) \to (t/\ell, x/\ell)$, which leaves the inviscid 
Burgers equation invariant, equation \eqref{h=1} becomes
\begin{equation}\label{ell=h=1}
u_t\ +\ u\,u_x\ =\ u_{xxt}\, -\,u_x\,u_{xx}\, +\,u\,u_{xxx}.
\end{equation}

Equation \eqref{ell=h=1} belongs to a three-parameter family 
of non-dispersive equations, given by 
\begin{equation}\label{generalequation}
u_t\ -\ u_{xxt}\ =\ a\, u\, u_x\ +\ b\, u_x\, u_{xx}\ +\ c\, u\, u_{xxx},
\qquad\mbox{for $a,b,c\in\mathds{R}$.} 
\end{equation}
In this family, we look for an equation that has Galilean invariance and 
conservation of energy (at least for smooth solutions).
The family \eqref{generalequation} includes a number of famous equations, including the 
the dispersionless Camassa--Holm equation \cite{CamassaHolm1993} 
\begin{equation}\label{CH}
u_t\ +\ 3\, u\, u_x\ =\  u_{xxt}\ +\ 2\, u_x\, u_{xx}\ +\  u\, u_{xxx},
\end{equation}
as well as the Degasperis--Procesi equation \cite{DegasperisProcesi1999} and the Benjamin--Bona--Mahony equation  
\cite{BenjaminEtAl1972}. 
None of these equations have the properties we seek, however. 
E.g., it is well known that the Camassa--Holm equation conserves the $H^1$ energy \cite{CamassaHolm1993}, 
but is not Galilean invariant. 

%
%

In order to obtain a Galilean invariant regularisation of the Burgers equation, one must take 
$c=-a=1$ in (\ref{generalequation}). The special case $b=0$ was studied by Bhat and Fetecau 
\cite{BhatFetecau2006,BhatFetecau2007,BhatFetecau2009a}, who proved the existence 
of the solution and the convergence to weak solutions of the Burgers equation
in the limit corresponding to $\ell\to0$. The limit fails 
to satisfy the entropy condition for the Riemann problem with $u_\text{left}\/ <\/ u_\text{right}$ 
\cite{BhatFetecau2009a}. For this regularisation, no energy conservation equation is known.

In the present paper, we consider $c=-a=1$ (to ensure Galilean invariance, as in 
\citep{BhatFetecau2006}) and, in order to maintain conservation of the $H^1$ norm at least for 
smooth solutions, we take $b=2$ (as in the Camassa--Holm equation). With this done, equation 
\eqref{generalequation} becomes 
\begin{equation}
u_t\ +\ u\, u_x\ =\ u_{txx}\ +\ 2\, u_x\ u_{xx}\ + u\, u_{xxx},
\end{equation}
Introducing the scaling $(t,x)\mapsto(\ell\, t, \ell\, x)$, we obtain 
\begin{equation}\label{rB}
u_t\ +\ u\, u_x\ =\ \ell^2\left( u_{txx}\ +\ 2\, u_x\, u_{xx}\ + u\, u_{xxx} \right),
\end{equation}
that is a formal approximation of the Burgers equation for small $\ell$. 
Equation (\ref{rB}) is the regularised Burgers (rB) equation studied in this paper.
It is Galilean invariant, and 
smooth solutions of (\ref{rB}) satisfy a conservation law for an $H^1$ energy density, namely
\begin{equation}\label{rBene0}
\left[\,\half\,u^2\,+\,\half\,\ell^2\,u_x^{\,2}\,\right]_t\ +\,\left[\,\third\,u^3\,
-\,\ell^2\,u^2\,u_{xx}\,-\,\ell^2\,u\,u_{xt}\,\right]_x\ =\ 0.
\end{equation}
A mathematically equivalent equation was proposed in \cite[Remark 1]{DDM}
as a modification of the BBM equation that possesses a Galilean-like invariance property.

We remark that the rB equation \eqref{rB} has variational structure (described further in \cite{Guelmame2020})
that we will not use here and which appears unrelated to $H^1$ energy conservation. 
E.g., \eqref{rB} can be obtained as the Euler--Lagrange equation for an action with Lagrangian density
\begin{align}\label{lagRB}
\mathcal{L}_\ell\ \eqdef\ \half\,\phi_x\,\phi_t\ +\ \sixth\,\phi_x^{\,3}\ +\ 
\half\,\ell^2\, \left[ \phi_x\,\phi_{xx}^{\,2}\ -\ \phi_{xxx}\,\phi_{t} \right],
\end{align}
where $\phi$ is a velocity potential, i.e., $u=\phi_x$.
The rB equation has also a Hamiltonian structure \cite{DDM,Guelmame2020}, but with Hamiltonian different from the $H^1$ energy.
Indeed, with the Hamiltonian operator and functional 
\begin{align}
\mathscr{D}\ \eqdef\ \left(1\/ -\/ \ell^2\, \partial_x^2\right)^{-1}\,\partial_x, \qquad
\mathfrak{H}\ \eqdef\, \int\left[\,\sixth\,u^3\,+\,\half\,\ell^2\,u\,u_x^{\,2}\,
\right]\,\ud\/x,
\end{align}
the equation
\begin{align}
u_t\ &=\ -\,\mathscr{D}\ \delta_u\,\mathfrak{H}\ =\ -\,\mathscr{D}\left[\,\half\,u^2\, 
+\,\half\,\ell^2\,u_x^{\,2}\,-\,\ell^2\,(u\/\/u_x\/)_x \,\right]
\end{align}
can be rewritten in a form equivalent to the rB equation \eqref{rB}, namely
\begin{align}\label{RB_2}
u_t\ +\ u\,u_x\ +\ \ell^2\,P_x\ =\ 0, \qquad
P\ \eqdef\ \mathfrak{G}\ast \half\, u_x^{\,2} \geqslant 0 , \qquad  \mathfrak{G}\ 
\eqdef\ (2\ell)^{-1}\,\exp(-|x|/\ell),
\end{align}
where $\ast$ denotes the convolution product. 
For comparison, the Camassa--Holm equation~\eqref{CH} can be rewritten 
in the form 
\begin{align}\label{CHconv}
u_t\ +\ u\,u_x\ +\, \left[\, \mathfrak{G} \ast \left( \half\, u_x^{\,2}\, +\ u^2 \right) 
\/\right]_x\  =\ 0, \qquad  \mathfrak{G}(x)\ \eqdef\ \half\,\exp(-|x|).
\end{align}

Differentiating \eqref{RB_2} with respect to $x$, and using that $P-\ell^2\/P_{xx}\/
=\/\half\/u_x^2$, one obtains
\begin{equation}\label{Ric}
\left[\, u_t\, +\, u\,u_x\, \right]_x\ +\ P\ = \half\, u_x^2.
\end{equation}
Note that $P$ goes formally to zero as $\ell \to +\infty$, whence one obtains the 
Hunter--Saxton (HS) equation \cite{HunterSaxton1991,HunterZheng1994}
\begin{equation}\label{HS}
\left[\, u_t\, +\, u\, u_x\, \right]_x\ =\ \half\, u_x^2.
\end{equation}
Note also that by taking $\ell \to +\infty$ formally in \eqref{rB}, we obtain the derivative 
of \eqref{HS} with respect of $x$.

In this section, we have provided a heuristic derivation of a regularised Burgers equation
by imposing the important physical requirements of Galilean invariance and energy conservation. 
We have also related this equation with well-known equations. 
In the rest of the paper, we perform a
rigorous mathematical investigation of its solutions.

\section{Existence and breakdown of smooth solutions}\label{Existence of the solution}

This section is devoted to show the local existence and breakdown of smooth solutions 
for the  Cauchy problem \eqref{RB_2} with $u(0,x)=u_0(x)$. The form \eqref{RB_2} of the 
regularised Burgers equation is more convenient for studying smooth solutions than (\ref{rB}), 
because it involves fewer derivatives. 

Usually, one needs an equation for $u_x$ to study the life span of smooth solutions.
Equation \eqref{Ric} can be written
\begin{equation}\label{RB_x}
u_{xt}\ +\ \half\, u_x^2\ +\ u\, u_{xx}\ +\ P\ =\ 0.
\end{equation}
Multiplying \eqref{RB_2} by $u$ and multiplying \eqref{RB_x} by $\ell^2\, u_x$, we obtain 
\begin{align}
\left[\,\half\,u^2\,\right]_t\ +\,\left[\,\third\,u^3\,+\,\ell^2\,u\,P\,\right]_x\ 
&=\  \ell^2\, u_x\, P,\label{rBene1} \\
\left[\,\half\,\ell^2\, u_x^2\,\right]_t\ +\, \left[\,\half\,\ell^2\,u\,u_x^2\,\right]_x\ 
&=\ -\,\ell^2\,u_x\,P,\label{rBene2}
\end{align}
which imply an energy conservation law for smooth solutions; i.e., we have 
the (conservative) energy equation
\begin{equation}\label{rBene}
\left[\,\half\, u^2\,+\,\half\,\ell^2\,u_x^2\,\right]_t\ +\, \left[\,\third\,u^3\,+\, 
\ell^2\,u\,P\,+\,\half\,\ell^2\,u\,u_x^2\,\right]_x\ =\ 0.
\end{equation}

For a class of equations including $\mathrm{rB}$ as special case, \citet{Yin2004b,yin2007cauchy} 
has proven the following local existence result.
\begin{thm}[\citet{Yin2004b,yin2007cauchy}]\label{Yin} 
For an initial datum $u_0 \in H^s(\mathds{R})$ with $s>3/2$, there exists a maximal time 
$T^*>0$ (independent of $s$) and a unique solution $u \in \mathcal{C}([0,T^*[,H^s)$ of 
\eqref{RB_2} such that (blow-up criterium)
\begin{equation}
T^* < +\infty \qquad \implies \qquad  \lim_{t \uparrow T^*}  \|\,u(t,\cdot)\,
\|_{H^s}\ =\ +\infty.
\end{equation}
Moreover, if $s \geqslant 3$, then
\begin{equation}
T^* < +\infty \qquad \implies \qquad   \lim_{t \uparrow T^*}\ \inf_{x\in\mathds{R}}\ u_x(t,x)\ =\ -\infty.
\end{equation}
\end{thm}
\noindent
Furthermore, the solution given in this theorem satisfies the Oleinik inequality:
\begin{pro}{\bf (Oleinik inequality)}\label{Oleiniki}
Let $u_0 \in H^s(\mathds{R})$ with $s\, \geqslant\,  2$ and let $M=\sup_{x \in 
\mathds{R}} u_0'(x)$. Then, for all $t \in [0,T^*[$ the solution given 
in Theorem \ref{Yin} satisfies
\begin{equation}\label{OlM}
\sup_{x\in \R} u_x(t,x)\ \leqslant\ { \textstyle \frac{2\,M}{M\,t\,+\,2}}\ \leqslant\ M.
\end{equation}
\end{pro}
\proof
Let $x_0 \in \mathds{R}$ and let the characteristic $\eta(t,x_0)$ be defined as the 
solution of the Cauchy problem $\eta_t(t,x_0)=u(t,\eta(t,x_0))$, with the initial 
datum $\eta(0,x_0)=x_0$. With $H(t,x_0)\eqdef u_x(t,\eta(t,x_0))$, the equation \eqref{RB_x} 
can be rewritten
\begin{equation}\label{Riccati}
H_t\ +\ \half\,H^2\ +\ P\ =\ 0.
\end{equation}
Since $P \geqslant 0$, it follows that  $H_t \leqslant -\half H^2$ which implies 
that $H(t,x_0) \leqslant \frac{2 H(0,x_0)}{H(0,x_0)t+2} \leqslant 
\frac{2M}{Mt+2}$. \qed
\begin{rem}
The Oleinik inequality \eqref{OlM} is valid only when the solution $u$ is smooth. 
In Theorem \ref{existence2} below, we show that this inequality holds 
for all times also for a certain type of weak solutions (called {\em dissipative})
such that $u \in H^1$ (and, possibly, for $M=+\infty$).
\end{rem}
Unfortunately, the solution given in Theorem \ref{Yin} does not exist globally 
in time for all non trivial initial data \citep{Yin2004b}. Since \citet{Yin2004b} 
studied a general family of equations including $\mathrm{rB}$, his result is not 
optimal for $\mathrm{rB}$. In the following proposition, this result is improved 
with a shorter proof. 
\begin{pro}{\bf (An upper bound on the blow-up time)}\label{blowup}
Let $u_0 \in H^s(\mathds{R})$ with $s\, \geqslant\, 2$. If there exists $x_0 \in 
\mathds{R}$ such that $u'_0(x_0) < 0$, then $T^* \leqslant -2/\inf u_0'.$
\end{pro}
\proof
From the proof of the previous proposition, we have 
\begin{equation}\label{H}
 H(t,x_0)\ \leqslant\ \frac{2 H(0,x_0)}{t H(0,x_0)+2}, 
 \qquad t<T^*.
\end{equation}
If $T^* > -2/\inf u_0'$ then $H(0,x_0)<0$ implies 
\[
\pushQED{\qed} 
\lim\limits_{t \to -2/H(0,x_0)} H(t,x_0)=-\infty,  
\] 
this contradicts 
$u \in \mathcal{C}([0,T^*[,H^s)$. \qed

A uniform (with respect to $\ell$) lower bound on $T^*$ 
is needed, in order to prove in section \ref{Compactness} below the convergence of 
smooth solutions (see Proposition \ref{mu} and Proposition \ref{nu}).
\begin{thm}{\bf (A lower bound on the blow-up time)}\label{InteractionTime}
Let $u_0$ in $H^s$ be non-trivial with $s\, \geqslant\, 2$ and let 
\begin{equation*}
m(t)\ \eqdef\  \inf_{x \in \mathds{R}} 
u_x(t,x)\ <\   0\ <\ M(t)\ \eqdef\ \sup_{x \in \mathds{R}} u_x(t,x), \qquad t\ <\ T^*.
\end{equation*}
If $|m(0)| \geqslant M(0)$ then
\begin{equation}
-1\left/\/\inf u_0'\right.\ \leqslant\ T^*\ \leqslant\ -2\left/\/\inf u_0'.\right.
\end{equation}
If $|m(0)| < M(0)$ then, there exists $t^*$ such that $0<t^*\leqslant-m(0)^{-1}-M(0)^{-1}$
and $m(t^*)=-M(t^*)$. Therefore
\begin{equation}
t^*\ +\ 1\left/\/\sup u_0'\right.\ \leqslant\ T^*\ \leqslant\ -2\left/\/\inf u_0'.\right.
\end{equation}
\end{thm}
\begin{rem}\label{BBT}
Note that the blow-up time $T^*$ is uniformly 
(with respect to $\ell$) bounded from below 
by $1/\sup |u_0'|$.
\end{rem}

\proof
Since $u\, \in\, H^s$, $u_x \rightarrow 0$ when $x$ goes to $\pm \infty$, and $u_x$ is not the zero function, 
so $m(t)= \min_{x \in \mathds{R}}H(t,x)<0<M(t)=\max_{x \in \mathds{R}}H(t,x)$.
The equation \eqref{Riccati} implies that $m$ and $M$ are decreasing in time, so $|m|\, =\, -m$ is increasing. 
So, if $|m(t_0)|\ \geqslant\ M(t_0)$, then for all $t>t_0$ we have $|m(t)|\ \geqslant\ M(t)$. 

The inequality \eqref{H} shows that for $t<T^*$
\begin{equation}\label{mM}
0\ <\ M(t)\ \leqslant\ \frac{2\, M(0)}{M(0)\, t\ +\ 2}, \qquad m(t)\ \leqslant\ \frac{2\, m(0)}{m(0)\, t\ +\ 2}\ <\ 0,
\end{equation}
which implies that, if $|m(0)| < M(0)$, there exists $t^*\leqslant-(m(0)+M(0))/(m(0)\/M(0))$ 
such that $|m(t^*)|=M(t^*)$.

If $\delta>0$ is small enough, since the function $H(t+\delta, \cdot)$ has a minimum, then there exists $x_\delta$ such that $m(t+\delta)=
H(t+\delta,x_\delta)$. Inspired by \citet{JuncaLombard19} one gets 
\begin{align}
m(t+\delta)\ &=\ H(t+\delta, x_\delta)\ =\ H(t,x_\delta)\ +\ \int_t^{t+\delta} 
H_t(s,x_\delta)\, \mathrm{d}s\nonumber \\
&\geqslant\ m(t)\ -\ \int_t^{t+\delta} \left( \half\, H(s,x_\delta)^2\ +\ P(s,x_\delta)  \right)\, \mathrm{d}s.
\end{align}
Since $m(\cdot) < 0$ and $\delta$ is arbitrary small, we have $m(s) \leqslant H(s,x_\delta) \leqslant 0$ then  $m(s)^2 \geqslant H(s,x_\delta)^2$, implying that
\begin{equation}
\frac{m(t+\delta)-m(t)}{\delta}\ \geqslant\ -\,\frac{1}{\delta} \int_t^{t+\delta} 
\left( \half\, m(s)^2  + \sup_{x \in \mathds{R}} P(s,x)  \right) \mathrm{d}s .
\end{equation}
Defining the generalised derivative 
\begin{equation}
\dot{m}(t)\ \eqdef\ \liminf_{\delta\to 0^+}\, \frac{m(t+\delta)\,-\,m(t)}{\delta},
\end{equation}
one can show that 
\begin{equation}\label{Re}
\dot{m}(t)\ +\ \half\, m(t)^2\ \geqslant\ - \sup_{x \in \mathds{R}} P(t,x).
\end{equation}
Using the definition of $P$ from \eqref{RB_2} and using that $\|\mathfrak{G}\|_1=1$, one obtains
\begin{equation}
\sup_{x \in \mathds{R}} P(t,x)\ \leqslant\ \half\, \|\,u_x(t,\cdot)\,\|_\infty^2\ \leqslant\ \max 
\left\{ \half\, M(t)^2, \half\, m(t)^2 
\right \}.
\end{equation} 
and the Riccati-like inequality \eqref{Re} becomes
\begin{equation}
\dot{m}(t)\ +\  m(t)^2\ \geqslant\ 0 \qquad t>t^*.
\end{equation}
Then $T^*-t^* \geqslant -1\, /\, m(t^*)\, =\, 1\, /\, M(t^*)$, and with \eqref{mM}, one obtains
\[
\pushQED{\qed} 
T^*\ \geqslant\  t^*\ +\ 1\left/\sup u_0'.\right.  \qedhere
\popQED
\]

\section{Global weak solutions: conservative case}\label{Global weak solutions}

Note that Proposition \ref{blowup} shows that, for $s\geqslant2$, we have 
$\lim_{t \uparrow T^*} \inf_{x\, \in\, \mathds{R}} u_x(t,x)=-\infty$ which 
implies that
\[ 
\lim_{t \uparrow T^*}\ \|\,u(t,\cdot)\,\|_{H^s}\ =\ +\infty . 
\]
Hence the space $H^s$ with $s\geqslant2$ is not the right space in order  
to obtain the global existence of the solution.

\citet{BressanConstantin2007a,BressanConstantin2007b} have proved the existence 
of two types of global solutions for the Camassa--Holm equation \eqref{CHconv} in $H^1$.
Using the formal energy equation \eqref{rBene}, a similar proof (of global 
existence of conservative and dissipative solutions in $H^1$) for $\mathrm{rB}$ 
can be done following \citep{BressanConstantin2007a,BressanConstantin2007b}. 
Another proof of existence of a dissipative solution, using the vanishing viscosity method,  
is given by \citet{ChenTian2009,XinZhang2000}.

In this paper, the existence theorem 
will be developed for solutions not vanishing as $|x| \to \infty$. Note that a major 
difference between the $\mathrm{rB}$ \eqref{RB_2} and the Camassa--Holm \eqref{CHconv} 
equations is that $u^2$ does not appear in the non-local term of $\mathrm{rB}$. This
allows us to get global existence for $\mathrm{rB}$ without asking $u$ to be in 
$L^2(\mathds{R})$. Moreover, in Theorem \ref{optimality} below, we show that asking 
$u_x \in L^2$ is optimal.

These remarks lead us to assert in the following the existence of two types
of solutions of $\mathrm{rB}$: conservative and dissipative. 
We start the analysis in this section by defining a conservative solution.
\begin{definition}\label{defcon}
A function $u$ is called a {\em conservative} solution of 
$\mathrm{rB}$ if
\begin{itemize}
\item The function $u$ belongs to $\mathrm{Lip}([0,T], L^2_{\rm loc})$ and $ u_x \in L^\infty ([0,T],L^2_{\rm loc})$ for all $T>0.$
\item $u$ satisfies the equation \eqref{RB_2}, with an 
initial data $u(0,x)=u_0(x)$.
\item $u$ satisfies \eqref{rBene} in the sense of distributions.
\end{itemize}
\end{definition}

\noindent It means that it is a weak solution conserving the energy, as smooth solutions.

\begin{rem}\label{energyeq}
The regularity $u_x \in L^\infty([0,T],L^2_{\rm loc})$ ensures that \eqref{rBene1} is satisfied. Thence, the equalities \eqref{rBene2} and \eqref{rBene} are equivalent.
\end{rem}

Introducing the homogeneous Sobolev space $\dot{H}^1(\mathds{R})=\{f:\|f'\|_2< +\infty \}$, 
we can state the theorem:
\begin{thm}\label{existence}
Let $u_0 \in \dot{H}^1(\mathds{R}) \cap L^\infty(\mathds{R})$. If there exists a Lipschitz 
function $\phi$ such that $\phi' \in L^1(\mathds{R})$ with  $u_0-\phi \in H^1(\mathds{R})$, 
then there exists a global conservative solution $u$ of \eqref{RB_2}, such that $u(t,\cdot)
-\phi \in H^1(\mathds{R})$ for all $t>0$. In addition, for all $T>0$
\begin{equation}\label{AntiOl}
\lim_{t \uparrow T}\, \inf_{x\in\mathds{R}}\, u_x(t,x)\ =\  -\infty 
\quad \implies \quad 
\lim_{t \downarrow T}\, \sup_{x\in\mathds{R}}\, u_x(t,x)\ =\ +\infty,
\end{equation}
and if $u_0 \in H^1$, then for almost all $t>0$
\begin{equation}\label{conservation}
\int_\mathds{R}\left[\,u(t,x)^2\,+\, \ell^2\ u_x(t,x)^2\,\right]\mathrm{d}x\ =\ 
\int_\mathds{R}\left[\, u_0(x)^2\, +\, \ell^2\ u_0'(x)^2\,\right] \mathrm{d}x.
\end{equation}
\end{thm}
\begin{rem}
This theorem covers also some solutions that do not have a limit when $|x| 
\to \infty$, such as $\phi(x)=u_0(x)=\cos\/\ln(x^2+1)$.
\end{rem}
\begin{rem}\label{remAntiOl}
Note that \eqref{AntiOl} implies that the Oleinik inequality \eqref{OlM} cannot hold 
after the appearance of singularities.
\end{rem}

\textit{Proof of Theorem \ref{existence}.}
In the special case $u_0 \in H^1(\mathds{R})$, the proof can be done following 
\citet{BressanConstantin2007a}. In the general case, the energy is modified as
\begin{equation}\label{E}
E(t)=\int_\mathds{R} \left[ u(t,x)-\phi(x) \right]^2+\ \ell^2\ 
u_x(t,x)^2\, \mathrm{d}x ,
\end{equation}
and the proof is done in several steps: 
\begin{itemize} 
\item In the first step we obtain an energy estimate of the solution in the Eulerian coordinates.
\item In the second step we define a mapping from the Eulerian to the Lagrangian coordinates where we obtain an equivalent semi-linear system of rB.
\item In the third and the fourth steps we prove the existence of global solutions of the equivalent system in the Lagrangian coordinates.
\item In step 5, we rewrite the solution of the equivalent system in the Eulerian coordinates and we show that it is a global conservative solution of rB.
\end{itemize}

{\bf Step 1: Formal energy estimate on the $x$-variable.}
Let $\tilde{u}(t,x)=u(t,x)-\phi(x) $. The equation \eqref{RB_2} can be rewritten
\begin{equation}\label{RB}
u_t\ +\ u\,u_x\ +\ \ell^2\,P_x\ =\ \tilde{u}_t\ +\ u\, u_x\ +\ \ell^2\,P_x\ =\ 0.
\end{equation}
Multiplying (\ref{RB}) by $\tilde{u}$, one gets
\begin{equation}\label{L2}
\left[\,\half\,\tilde{u}^2\,\right]_t\ +\,\left[\,\third\,u^3\,-\,\half\,\phi\,u^2\, 
\right]_x\ +\ \half\,\phi_x\,u^2\ +\ \ell^2\,u\,P_x\ -\ \ell^2\,\phi\,P_x\ =\ 0.
\end{equation}
Adding \eqref{L2} and \eqref{rBene2}, we obtain
\begin{equation}
\half \left[\/\tilde{u}^2\/+\/\ell^2\,u_x^2\/\right]_t\, +\left[\/\third\/ 
u^3\/+\/\sixth\/\phi^3\/-\/\half\/\phi\/u^2\/+\/\half\/\ell^2\/u\/u_x^{\,2}\/+\/ 
\ell^2\/u\/P\/\right]_x\, =\, \ell^2\/\phi\/P_x\, -\, \half\/\phi_x \left(\/\tilde{u}^2
\/+\/2\/\phi\/\tilde{u}\/\right).
\end{equation}
Integrating over the real line, one gets (exploiting the triangular inequality)
\begin{equation}\label{preGronwall}
\half\,E'(t)\ \leqslant\  \int_\mathds{R} \left(\,\ell^2\,|\phi\/P_x|\,+\, 
\half\,|\phi'| \left(\,2\,\tilde{u}^2\,+\,\phi^2\,\right)\/  \right)\mathrm{d}x .
\end{equation}
The Young inequality implies that
\begin{subequations}\label{P,P_x}
\begin{align}
\|\,P(t)\,\|_p\  & \leqslant\ {\textstyle \frac{1}{2\,\ell^2}}\,\|\,\mathfrak{G}\,\|_p\,E(t) \qquad 
\forall p \in [1,\infty], \label{Young}\\
\|\,P_x(t) \,\|_p\ & \leqslant\ {\textstyle \frac{1}{2\,\ell^3}}\, 
\|\,\mathfrak{G}\,\|_p\,E(t) \qquad \forall p \in [1,\infty].\label{Young_x}
\end{align}
\end{subequations}
Using \eqref{preGronwall} and \eqref{Young_x}, we obtain
\begin{equation}
E'(t)\ \leqslant\,\left(\,\ell^{-1}\,\|\,\phi\,\|_\infty\, +\, 2\,\|\,\phi'\,\|_\infty\, \right) 
E(t)\ +\ \|\,\phi\,\|_\infty^2 \,\|\,\phi'\,\|_1.
\end{equation}
Then the Gronwall lemma ensures that $E(t)$ does not blow up in finite 
time. \\

{\bf Step 2: Equivalent system.}
As in \citep{BressanConstantin2007a}, let $\xi \in \mathds{R}$ and let $y_0(\xi)$ be defined by
\begin{equation}\label{y}
\int_0^{y_0(\xi)}\left(\,1\,+\,{u_0^\prime}^2\,\right) \mathrm{d}x\ =\ \xi,
\end{equation}
and let $y(t,\xi)$ be the function\footnote{It will turn out that $y(t,\xi)$ is 
the characteristic of $\mathrm{rB}$ corresponding to $y_0(\xi)$, with speed $u(t,y(t,\xi))$.} 
defined by the equation
\begin{equation}\label{char}
y_t(t,\xi)\ =\ u(t,y(t,\xi)), \qquad y(0,\xi)\ =\ y_0(\xi).
\end{equation}
Let also $v=v(t,\xi)$ and $q=q(t,\xi)$ be defined as
\begin{equation}\label{vq}
v\ \eqdef\ 2\, \arctan(u_x), \qquad q\ \eqdef\,\left(\,1\,+\,u_x^{\,2}\,\right)y_\xi,
\end{equation}
where $u_x(t,\xi) =  u_x(t,y(t,\xi))$. Notice that
\begin{equation}\label{identities}
\frac{1}{1+u_x^{\,2}}\ =\ \cos^2\!\left(\frac{v}{2}\right), \quad \frac{u_x}{1+u_x^{\,2}}
\ =\ \frac{\sin(v)}{2}, \quad \frac{u_x^{\,2}}{1+u_x^{\,2}}\ =\ \sin^2\!\left(\frac{v}{2}
\right), \quad 
\frac{\partial y}{\partial \xi}\ =\ q\/\cos^2\!\left(\frac{v}{2}\right).
\end{equation}
Integrating the last equality in \eqref{identities}, one obtains
\begin{equation}\label{cov2} 
y(t,\xi')\ -\ y(t,\xi)\ =\ \int_\xi^{\xi'} q(t,s)\/ \cos^2\!\left(\frac{v(t,s)}{2}\right) 
\mathrm{d}s.
\end{equation}
Using \eqref{identities} and the change of variables $x=y(t,\xi')$, \eqref{cov2}, $P$ and 
$P_x$ can be written in the new variables as
\begin{align} \nonumber
P(t,\xi)\ &=\ {\frac{1}{4\,\ell}} \int_\mathds{R} 
\exp\!\left(\frac{-|y(t,\xi)-x|}{\ell}\right) u_x^{\,2}(t,x)\,\mathrm{d}x  \\ \label{P} 
&=\ {\frac{1}{4\, \ell}} \int_\mathds{R} \exp\!\left(- \frac{1}{\ell} \left| \displaystyle 
\int_\xi^{\xi'} q(t,s) \cos^2\!\left(\frac{v(t,s)}{2}\right) \mathrm{d}s \right| \right) 
q(t,\xi') \sin^2\!\left(\frac{v(t,\xi')}{2}\right) \mathrm{d}\xi', \\ \nonumber
P_x(t,\xi)\ &=\ {\frac{1}{4\,\ell^2}} \left(\int_{y(t,\xi)}^{+\infty}-\int_{-\infty}^{y(t,\xi)} 
\right)\exp\!\left({-|y(t,\xi)-x|}{\ell}\right)  u_x^{\,2}(t,x)\,\mathrm{d}x \\  \label{P_x} 
&=\, \left( 
\int_{\xi}^{+\infty}-\int_{-\infty}^{\xi} \right)\ \exp\! 
\left(-\left| \displaystyle \int_\xi^{\xi'} q(t,s) \cos^2\!\left( 
\frac{v(t,s)}{2}\right) \frac{\mathrm{d}s}{\ell} \right| \right) q(t,\xi') \sin^2\!\left( 
\frac{v(t,\xi')}{2}\right) \frac{\mathrm{d} \xi'}{4\,\ell^2}.
\end{align}
Then, a system equivalent to the $\mathrm{rB}$ equation is given by
\begin{subequations}\label{ODE}
\begin{align}
y_t\ &=\ u,  &y(0,\xi)\ =\ y_0(\xi), \label{ODEa}\\ 
u_t\ &=\ - \ell^2\, P_x,  &u(0,\xi)\ =\ u_0(y_0(\xi)), \label{ODEb}\\ 
v_t\ &=\ -P\left(1+\cos(v)\right)\ -\ \sin^2({v}/{2}),  &v(0,\xi)\ =\ 2\/\arctan\!\left(
u_0' \left(y_0(\xi)\right)\right),\label{ODEc}\\ \
q_t\ &=\ q\left(\half-P\right) \sin(v),   & q(0,\xi)\ =\ 1. \label{ODEd}
\end{align}
\end{subequations}
In order to prove Theorem \ref{existence}, we prove first the global existence 
of the solution of the initial-value problem \eqref{ODE}, 
then we infer that this solution yields a conservative 
solution of $\mathrm{rB}$.\\

{\bf Step 3: Local existence for the new system.}
Our goal is to prove that the system of equations \eqref{ODE} is locally 
well-posed. The proof given in \citep{BressanConstantin2007a} for the 
Camassa--Holm equation is slightly simplified here. 

We first solve a coupled 2x2 subsystem instead of a 3x3 subsystem in \citep{BressanConstantin2007a}.
Let $u_0$ be a function such that $u_0-\phi \in H^1$, then $y_0$ is well 
defined in \eqref{y}. Note that the right-hand side of \eqref{ODE} does 
not depend on $y$. Since $P$ and $P_x$ depend only on $v$ and $q$, the 
right-hand sides of equations \eqref{ODEb}, \eqref{ODEc} and \eqref{ODEd} 
do not depend on $u$. Also, the equations \eqref{ODEc} and \eqref{ODEd} 
are coupled. Thus, we are left to show that the system of two equations 
\begin{subequations}\label{ODE2}
\begin{align}
v_t&=-P(1+\cos v)-\sin^2 \frac{v}{2},  &v(0,\xi) =v_0(\xi)=2\ \arctan 
u_0' \left(y_0(\xi)\right),\\
q_t&=q\ (\half-P) \sin v,   &q(0,\xi)=q_0(\xi)=1
\end{align}
\end{subequations}
is well defined in the space $X\eqdef\mathcal{C}([0,T],L^\infty(\mathds{R},
\mathds{R}^2))$.

Let $U=(v,q)$, and let $\mathfrak{D} \subset X$ be the closed set 
satisfying $U(0,\xi)=U_0(\xi)$ and
\begin{subequations} \label{eqUDC}
\begin{align}
1\, /\, C\ \leqslant\ q(t,\xi)\ &\leqslant\ C\ \quad \forall (t,\xi) \in 
[0,T] \times \mathds{R}, \\
\left| \left\{ \xi,\ \sin^2 {\textstyle  \frac{v(t,\xi)}{2}} \geqslant 
{\textstyle \half}  \right\} \right|\ & \leqslant\ C\ \quad \forall t \in 
[0,T],
\end{align}
\end{subequations}
where $C>0$ is a constant. Then, for $\xi_1<\xi_2$, from the equations 
(\ref{eqUDC}) we obtain
\begin{equation}
\int_{\xi_1}^{\xi_2} q(\xi)\, \cos^2 \frac{v(\xi)}{2}\, \mathrm{d} \xi\
\geqslant\ \int_{\left\{ \xi \in [\xi_1,\xi_2],\ \sin^2 {\textstyle  \frac{v(t,\xi)}{2}} \
\leqslant\  {\textstyle \half}  \right\}} \frac{C^{-1}}{2}\, \mathrm{d}\/\xi
\ \geqslant\, \left[ \frac{\xi_2-\xi_1}{2}-\frac{C}{2}  \right] C^{-1}.
\end{equation}
Let $\Gamma$ be defined as 
\begin{equation}\label{Gamma}
\Gamma(\zeta)=\min \left\{1,\ \exp \left( \frac{1}{2\, \ell}  
-\frac{|\zeta|}{2\, \ell} C^{-1} \right) \right\}.
\end{equation}
Then, for $(v,q) \in \mathfrak{D}$, the exponential terms in \eqref{P} and 
\eqref{P_x} are smaller than $\Gamma(\xi-\xi')$.

Let $P(\xi,v,q)$ be defined by \eqref{P}. If $(v,q) \in \mathfrak{D}$ then, 
using Young inequality, $\partial_v P$ and $\partial_q P$ are bounded, i.e., 
for $\{U, \tilde{U}\} \in \mathfrak{D}$ we have
\begin{equation}\label{PLip}
\|P(\xi,U)\ -\ P(\xi,\tilde{U})\, \|_X\ \lesssim\ \|U\ -\ \tilde{U}\,\|_X,
\end{equation}
where the symbol $\lesssim $ means ``less or equal'' with a constant  
depending only on $C$ and $\ell$.
Then, for $T$ small enough, the Picard operator
\begin{equation}
(\mathcal{P}(U))(t,\xi)\ \eqdef\ U_0\ +\ \int_0^t \left(- (1 + \cos v)\, P\ -\ \sin^2 
\frac{v}{2}\ ,\ q\, (\half-P)\, \sin v \right)\, \mathrm{d} \tau,
\end{equation}
is a contraction from $\mathfrak{D}$ to $\mathfrak{D}$. The local 
existence of the solution of the Cauchy problem \eqref{ODE2} follows at once.\\

{\bf Step 4: Global existence for the equivalent system.}
After proving the local existence of the solution of system 
\eqref{ODE2}, an estimate of the quantity
\begin{equation}\label{global}
\|\, q(t)\, \|_\infty\ +\ \left\|\, 1\, /\, q(t)\, \right\|_\infty\ +\ \left\|\, \sin^2  \left(
 v(t)\, / \, 2 \right) \, \right\|_1\ +\ \|\, v(t)\, \|_\infty ,
\end{equation}
is needed to ensure that the solutions exist for all time.
Let $u$ be defined as
\begin{equation}\label{u2}
u(t,\xi)\ \eqdef\ u_0(y_0(\xi))\ -\ \int_0^t \ell^2\, P_x(s,\xi)\, \mathrm{d}s,
\end{equation}
and let $y$ be  the family of characteristics
\begin{equation}\label{char2}
y(t,\xi)\ \eqdef\ y_0(\xi)\ +\ \int_0^t u(s,\xi)\, \mathrm{d} s,
\end{equation}
and, finally, let $\phi(t,\xi) \eqdef \phi(y(t,\xi))$. Our task here is to show 
that the modified energy
\begin{equation}\label{tilde{E}}
\tilde{E}(t)\ =\ \int_\mathds{R} \left[ (u-\phi)^2\, \cos^2 \frac{v}{2}\ +\  
\ell^2\, \sin^2 \frac{v}{2}\, \right]\, q\, \mathrm{d} \xi
\end{equation}
does not blow-up in finite time.

The system \eqref{ODE} implies that
\begin{equation}\label{E1}
\left( q\, \cos^2 \frac{v}{2} \right)_t\, =\ \half\, q\, \sin v, \qquad
\left( q\, \sin^2 \frac{v}{2} \right)_t\, =\ q_t\ -\ \half\, q\, \sin v\ =\ -\, q\,  P\, \sin v,
\end{equation}
while the equations \eqref{P} and \eqref{P_x} imply that
\begin{equation}\label{P_xi}
P_\xi\ =\ q\,  P_x\, \cos^2 \frac{v}{2}, \qquad \ell^2\, (P_x)_\xi\ =\ q\, P\, \cos^2 
\frac{v}{2}\ -\ \half\, q\, \sin^2 \frac{v}{2}.
\end{equation}
From \eqref{ODE}, \eqref{u2} and \eqref{P_xi}, we have
\begin{equation*}
\left( u_\xi\ -\ \half\, q\, \sin v \right)_t\ =\ 0,
\end{equation*}
and, for $t=0$, we have from \eqref{vq} and \eqref{identities}
\begin{equation*}
u_\xi\ -\ \half\, q\, \sin v\ =\ u_x\, {\textstyle{\partial y \over \partial 
\xi}}\ -\ \half\, \sin v\ =\ 0.
\end{equation*}
Thus, as long as the solution of \eqref{ODE} is defined, the equality
\begin{equation}\label{u_xi}
u_\xi\ =\ \half\, q\, \sin v
\end{equation}
holds. Therefore, the equations \eqref{E1}, \eqref{P_xi}, \eqref{u_xi} yield 
\begin{equation}\label{Ee1} 
\left[ \left( u^2\, \cos^2 \frac{v}{2}\ +\ \ell^2\, \sin^2 \frac{v}{2} 
\right)\, q\, \right]_t\, +\ \left[ 2\, \ell^2\, u\, P\ -\ \third\, u^3\, \right]_\xi\, =\ 0,
\end{equation}
which expresses conservation of energy in the $(t,\xi)$-variables when 
$u_+=u_-=0$, i.e., for $\phi=0.$

From \eqref{E1}, \eqref{u_xi} and \eqref{char2}, we have
\begin{equation*}
\left( q\, \cos^2 \frac{v}{2}\, \right)_t\ =\  u_\xi\ =\ \left( \frac{\partial 
y}{\partial \xi} \right)_t,
\end{equation*}
implying that the equality
\begin{equation}\label{y_xi}
\frac{\partial y}{\partial \xi}\ =\ q\, \cos^2 \frac{v}{2}
\end{equation}
holds for the $(t,\xi)$-variables (note that the equality is true for 
$t=0$ from \eqref{identities}). Then, using \eqref{char2} and \eqref{y_xi}, 
we get
\begin{equation}\label{phi}
\phi_t(t,\xi)\ =\ \frac{\ud}{\ud\/t}\, \phi(y(t,\xi))\ =\ u\, \phi', \qquad 
\phi_\xi(t,\xi)\ =\ q\, \cos^2\!\left(\frac{v}{2}\right) \phi',
\end{equation}
so, using \eqref{u2}, \eqref{ODE}, \eqref{phi} and \eqref{E1}, we obtain
\begin{align} \nonumber
\left[ (\phi^2\, -\, 2\, u\, \phi)\, q\, \cos^2 \frac{v}{2} \right]_t\, &+\, \left[ \half\, \phi\, 
u^2\, +\, \half\, \phi\, (u\, -\, \phi)^2\, -\, \sixth\, \phi^3\, \right]_\xi\ = \\ \nonumber
&2\, \ell^2\, P_x\, \phi\, q\, \cos^2 \frac{v}{2}\, -\, 2\, u\, \phi_\xi\, (u\, -\, \phi)\, +\, \half\, \phi_\xi\, u^2\, +\, (\phi^2\, -\, 2\, u\, \phi)u_\xi\, \\ \nonumber  
&+\, \phi\, u\, u_\xi\ +\, \half\, \phi_\xi\, (u\, -\, \phi)^2\, +\, \phi\, (u\, -\,\phi)_\xi (u\, -\, \phi)\
-\, \half\,  \phi^2 \phi_\xi  \\  \label{Ee2}
=\, &-\, \phi_\xi (u\, -\, \phi)^2+\, 2\, \ell^2\, P_x\, \phi\, q\, \cos^2 \frac{v}{2}\, -\, 2\, \phi\, \phi_\xi(u\, -\, \phi).
\end{align}
Adding \eqref{Ee1} and \eqref{Ee2}, with the trivial relation $2\, \phi\, (u\ -\ \phi)\ 
\leqslant\  \phi^2\ +\ (u\ -\ \phi)^2$, then integrating the result with respect of $\xi$, 
we get
\begin{equation}
\tilde{E}'(t)\ \leqslant\ \int_\mathds{R} \left( 2\, \ell^2\, |\phi\, P_x|\, q\, 
\cos^2 \frac{v}{2}\ +\  |\phi_\xi|\, \left( 2\, (u\ -\ \phi)^2\, +\ \phi^2 \right) 
\right)\, \mathrm{d} \xi .
\end{equation} 

Using \eqref{y_xi} and \eqref{phi} with the change of variables $x=y(t,\xi)$, 
then expoiting \eqref{P,P_x}, one obtains 
\begin{align*}
\tilde{E}'(t)\ &\leqslant\  \int_{ \{ \xi, \cos v \neq -1 \} } \left( 2\, \ell^2\, |\phi\, P_x|\ +\  |\phi'|\, \left( 2\, (u\ -\ \phi)^2\, +\ \phi^2 \right) 
\right)\, q\, 
\cos^2 \frac{v}{2}\, \mathrm{d} \xi \\
&=\ \int_\mathds{R} \left( 2\, \ell^2\, |\phi\, P_x|\ +\  |\phi'|\, \left( 2\, (u\ -\ \phi)^2\, +\ \phi^2 \right) 
\right)\, \mathrm{d} x \\
&\leqslant\ \left(\,\ell^{-1}\,\|\,\phi\,\|_\infty\, +\, 2\,\|\,\phi'\,\|_\infty\, \right) 
E(t)\ +\ \|\,\phi\,\|_\infty^2 \,\|\,\phi'\,\|_1,
\end{align*}
where $P_x$ in the second equation is defined as $P_x\ =\ \half\, \mathfrak{G}_x 
\ast u_x^2$.

From \eqref{E} and \eqref{tilde{E}}, and using the change of variables $x=y(t,\xi)$, 
one can show easily that 
\begin{equation}
E(t)\ =\ \int_{ \{ \xi, \cos v \neq -1 \}} \left[ (u-\phi)^2\, \cos^2 \frac{v}{2}\ +\  
\ell^2\, \sin^2 \frac{v}{2}\, \right]\, q\, \mathrm{d} \xi\ \leqslant\ \tilde{E}(t).
\end{equation}
Thence, the uniform estimate of $\tilde{E}(t)$ on 
any bounded interval $[0,T]$ follows by using Gronwall lemma.

We can show now that the quantity \eqref{global} does not blow up in finite time.
Using Young inequality, \eqref{P}, \eqref{P_x} and \eqref{Gamma}, one obtains 
\begin{subequations}\label{P,P_x2}
\begin{align}
\| P(t) \|_p\  & \leqslant\ \textstyle \frac{1}{4\, \ell^3}\, \|\Gamma\|_p\, \tilde{E}(t) 
\qquad \forall p \in [1,\infty], \\
\| P_x(t) \|_p\ & \leqslant\ \textstyle \frac{1}{4\, \ell^4}\, \|\Gamma\|_p\, 
  \tilde{E}(t) \qquad \forall p \in [1,\infty].
\end{align}
\end{subequations}
The inequalities \eqref{P,P_x2} are the identical estimates as \eqref{P,P_x}, but 
in the $(t,\xi)$-variables. Using \eqref{ODEd} and \eqref{P,P_x2}, we get
\begin{equation}
|q_t|\ \leqslant\ \left(\half\ +\ \textstyle \frac{1}{4\, \ell^3}\, E(t) \right)\, q,
\end{equation}
implying that $\|\/q(t)\/\|_\infty\/ +\/ \left\|\/ 1\, /\, q(t)\/ \right\|_\infty$ 
does not blow-up in finite time.
The equation \eqref{ODEc} and \eqref{P,P_x2} imply that 
$\|\/ v(t)\/ \|_\infty$ remains bounded on any finite interval $[0,T]$.
Also, the boundedness of the energy $\tilde{E}(t)$ and
$\left\|\, 1\, /\, q(t)\, \right\|_\infty$ implies that
$\left\|\, \cos^2 \left( v(t)\, /\, 2 \right) \,  \right\|_1$ remains bounded 
on any interval $[0,T]$. This completes the proof of the global existence. \\

{\bf Step 5: Global existence of a conservative solution.}
Here, we show that the global solution of the equivalent system (\ref{ODE}) 
yields a global solution of the $\mathrm{rB}$ equation.

Let $u$ and $y$ be defined by \eqref{u2} and \eqref{char2}, respectively. 
We claim that the solution of $\mathrm{rB}$  can be written as
\begin{equation}\label{solution}
u(t,x)\ =\ u(t,\xi), \qquad y(t,\xi)\ =\ x.
\end{equation}
Using \eqref{u_xi}, \eqref{phi} and the change of variables $x=y(t,\xi)$ 
with \eqref{y_xi}, one obtains
\begin{align*}
|u(t,\xi)\ -\ \phi(t,\xi)|^2\ &\leqslant\ 2\, \int_\mathds{R} |u\ -\ \phi|\, 
|u_\xi\ -\ \phi_\xi|\, \mathrm{d} \xi \\
&\leqslant\  2\, \int_\mathds{R} |u\ -\ \phi|\, q\, \left( \sin \frac{v}{2}\, \cos 
\frac{v}{2}\ +\ \phi'\, \cos^2 \frac{v}{2} \right)\, \mathrm{d} \xi \\
&\leqslant 2\, E(t)\ +\ \|\phi'\|_2^2,
\end{align*}
implying that $\|u(t)\|_\infty$ is uniformly bounded on any bounded 
interval $[0,T]$. Therefore, from \eqref{char2}, we get
\begin{equation}
y_0(\xi)\ -\ \|u(t)\|_\infty\, t\ \leqslant\  y(t,\xi)\ \leqslant\ 
y_0(\xi)\ +\ \|u(t)\|_\infty\, t,
\end{equation}
and thus
\begin{equation}
\lim\limits_{\xi \to \pm \infty} y_0(t,\xi)\ =\ \pm \infty .
\end{equation}
The equation \eqref{y_xi} implies that the mapping $\xi \mapsto y(t,\xi)$ is 
non-decreasing and, if for $\xi\ <\ \xi'$ we have $y(t,\xi)\ =\ y(t,\xi')$, 
then $\sin(v) = 2 \cos (v/2) \sin (v/2) = 0 $ between $\xi$ and $\xi'$ (see 
eq. \ref{y_xi}).
Integrating \eqref{u_xi} with respect to $\xi$, one obtains that 
$u(t,\xi) = u(t,\xi')$, so $u$ is well-defined in \eqref{solution}.

Proceeding as in \cite[section 4]{BressanConstantin2007a}, we can prove that 
for each interval $[t_1,t_2]$ there exists a constant $C = C(\ell,t_2)$  
such that, $\forall t \in [t_1\, ,t_2\ -\ h]$,
\begin{equation}\label{Lip}
\int_\mathds{R} |u(t\, +\, h,x)\ -\ u(t,x)|^2 \mathrm{d}x\ \leqslant\ C\, h^2 ,
\end{equation}
and then $u$ satisfies
\begin{equation}
\frac{\mathrm{d}}{\mathrm{d}t}\, u(t,y(t,\xi))\ =\ -\, P_x(t,\xi).
\end{equation}
{The inequality \eqref{Lip} implies that $u$ belongs to $Lip([0,T],L^2_{loc}).$}
Straightforward calculations show that, for $x=y(t,\xi)$ and for $\cos (v(t,\xi)) \neq 1$, 
we have 
\begin{equation}\label{u_x}
u_x(t,x)\ =\ \tan\!\left(\frac{v(t,\xi)}{2}\right)\, =\ \frac{\sin 
(v(t,\xi))}{1\,+\,\cos(v(t,\xi))}.
\end{equation}
Using the change of variables $x=y(t,\xi)$ with \eqref{y_xi}, one can 
show that $u$ is a global solution of $\mathrm{rB}$.

In order to prove \eqref{rBene2}, let $\psi$ be a test function and let 
$\tilde{\psi}(t,\xi) = \psi(t,y(t,\xi))$. Multiplying (\ref{E1}{\it b}) by $\tilde{\psi}$ 
and integrating the result with respect to $\xi$, one obtains
\begin{align} \nonumber
0\ =&\ \int_0^{+\infty} \int_{\mathds{R}} \left[ \left( q\, \sin^2 v/2 \right)_t\, +\, q\, P\, \sin v \right] \tilde{\psi}\, \mathrm{d}t\, \mathrm{d} \xi, \\ \nonumber
=&\ \int_0^{+\infty} \int_{\mathds{R}} \left[ -\tilde{\psi}_t\,  q\, \sin^2 v/2\, +\, \tilde{\psi}\, q\, P\, \sin v \right] \mathrm{d}t\, \mathrm{d} \xi\ +\ \int_{\mathds{R}} \tilde{\psi}(0,x)\, \sin^2 v(0,\xi)/2\, \mathrm{d}\xi, \\ \nonumber
=&\ \int\! \int_{\{\cos v > -1 \}} \left[ -\tilde{\psi}_t\,  q\, \sin^2 v/2\, +\, \tilde{\psi}\, q\, P\, \sin v \right] \mathrm{d}t\, \mathrm{d} \xi\ 
+\ \int_{\{v_0 > -\pi\}} \tilde{\psi}(0,x)\, \sin^2 v(0,\xi)/2\, \mathrm{d}\xi, \\ \label{localene}
&+\ \int\! \int_{\{\cos v = -1\}} -q\, \tilde{\psi}_t\, \mathrm{d}t\, \mathrm{d} \xi\ 
+\ \int_{\{v_0 = -\pi\}} \tilde{\psi}(0,x)\, \mathrm{d}\xi.
\end{align}
It is clear from \eqref{ODEc} that
\begin{equation}\label{measure}
\left| \left\{ \xi, \cos v(t,\xi)=-1 \right\} \right|=0 \quad 
\textrm{for almost all } t \geqslant 0.
\end{equation}
Using that $\tilde{\psi}_t = \psi_t + u\/\psi_x$ and the change of variables $x=y(t,\xi)$, 
the equation \eqref{rBene2} follows in the sense of distributions.

Finally, let $u_0 \in H^1$. The equation \eqref{Ee1} implies that
\begin{equation}\label{conservation2}
\frac{\mathrm{d}}{\mathrm{d}t} \int_\mathds{R} \left(  u^2\, \cos^2 
{\textstyle \frac{v(t,\xi)}{2}}\ +\ \ell^2\, \sin^2 {\textstyle 
\frac{v(t,\xi)}{2}} \right)\, q(t,\xi)\,  \mathrm{d}  \xi\ =\ 0,
\end{equation}
hence $\tilde{E}(t)=\tilde{E}(0)$. In addition, using the change of variables 
$x=y(t,\xi)$ with \eqref{y_xi} and \eqref{u_x}, one obtains
\begin{equation}\label{energyrelation}
\int_\mathds{R} u(t,x)^2\ +\ \ell^2\, u_x(t,x)^2\, \mathrm{d}x\
=\, \int_{\{ \xi,\ \cos v(t,\xi)>-1  \}} \left(  u^2\, \cos^2 {\textstyle 
\frac{v(t,\xi)}{2}}\ +\ \ell^2\, \sin^2 {\textstyle \frac{v(t,\xi)}{2}} 
\right)\, q(t,\xi)\,  \mathrm{d}  \xi .
\end{equation}
Using \eqref{measure}, the conservation of the energy \eqref{conservation} follows. 

We end this demonstration with the proof of the property \eqref{AntiOl}. The equation 
\eqref{ODEc} implies that $v$ is decreasing in time. Further, if $v(T,\xi)=-\pi$ (corresponding 
to an infinite value of $u_x$, see \eqref{u_x} above) then $v_t(T,\xi)=-1$, meaning that 
the value of  $v(t,\xi)$ crosses $-\pi$ and $v(t,\xi) < -\pi$ for all $t>T$. Then, 
\eqref{AntiOl} follows using \eqref{u_x}.\qed \\

\section{Global weak solutions: dissipative case}\label{ssecglobexidis}

We start this section by defining {\em dissipative\/} solutions,   
this kind of solution being very important for applications. 
We note in passing that when $\ell$ goes to zero, we expect to recover the entropy 
solution of the Burgers equation. 
{However, in Section \ref{Compactness}, we show that the limit (up to a subsequence) 
is a solution of the Burgers equation with a remaining forcing term.}

\begin{definition}\label{defdis}
A function $u$ is called a {\em dissipative} solution of 
$\mathrm{rB}$ if
\begin{itemize}
\item The function $u$ belongs to $\mathrm{Lip}([0,T], L^2_{\rm loc})$ and $ u_x \in L^\infty ([0,T],L^2_{\rm loc})$ for all $T>0.$;
\item $u$ satisfies the equation \eqref{RB_2}, with an 
initial data $u(0,x)=u_0(x)$;
\item $u$ satisfies the inequality 
\begin{equation}\label{localdiss}
\left[\,\half\, u^2\,+\,\half\,\ell^2\,u_x^2\,\right]_t\ +\, \left[\,\third\,u^3\,+\, 
\ell^2\,u\,P\,+\,\half\,\ell^2\,u\,u_x^2\,\right]_x\ \leqslant\ 0,
\end{equation}
in the sense of distributions.
\item {There exists a constant $C$ such that $u$ satisfies the Oleinik inequality 
\[
u_x(t,x)\ \leqslant\ C\,/\,t \qquad \forall t, x.
\]}
\end{itemize}
\end{definition}
\begin{rem}
{Following \citep{BressanConstantin2007b}, we construct in Theorem \ref{existence2} a dissipative 
solution of rB with $C=2$. The entropy solutions of the classical Burgers equation satisfy the 
Oleinik inequality with $C=1$.}
\end{rem}

As mentioned above, when $v$ crosses the value $-\pi$, $u_x$ jumps from $-\infty$ to $+\infty$, 
which means that the Oleinik inequality cannot be satisfied. Thus, to enforce the Oleinik inequality, 
the value of $v$ is not allowed to leave the interval $[-\pi,\pi[$. For that purpose, the system \eqref{ODE2} 
is modified (as in \citep{BressanConstantin2007b}) to become 
\begin{subequations}\label{ODE3}
\begin{align} \label{ODE3a}
u_t\ &=\ -\ell^2\, P_x,\\ \label{ODE3b}
v_t\ &=\ \begin{cases}
       -P\, (1+\cos v)\ -\ \sin^2(v/2), & v\,>\,-\pi, \\
       0 & v\, \leqslant\, -\pi ,
    \end{cases}\\ \label{ODE3c}
q_t\ &=\ \begin{cases}
      q\, (\half-P)\, \sin(v),   & v\,>\,-\pi \\
       0 & v\, \leqslant\, -\pi.
    \end{cases}
\end{align}
\end{subequations}
and $P$ and $P_x$ are also modified as  {\small
\begin{align}\label{P_dis}
P(t,\xi)\ &=\ { \frac{1}{4\, \ell} } \int_\mathds{R} \exp 
\left\{- \frac{1}{\ell}\, \Big| \displaystyle \int_\xi^{\xi'} \bar{q}(t,s)\,
\cos^2 \frac{v(t,s)}{2}\, \mathrm{d}s \Big| \right\}\,  \bar{q}(t,\xi')\, 
\sin^2 \frac{v(t,\xi')}{2}\, \mathrm{d} \xi', \\ \label{P_x_dis}
P_x(t,\xi)\ 
&=\ { \frac{1}{4\, \ell^2} } \left( 
\int_{\xi}^{+\infty}-\int_{-\infty}^{\xi} \right)\ \exp 
\left\{-\frac{1}{\ell}\, \Big| \displaystyle \int_\xi^{\xi'} \bar{q}(t,s)\,
\cos^2 \frac{v(t,s)}{2}\, \mathrm{d}s \Big| \right\}\, \bar{q}(t,\xi')\,
\sin^2 \frac{v(t,\xi')}{2}\, \mathrm{d} \xi',
\end{align}}
where $\bar{q}(t,\xi) = q(t,\xi)$ if $v(t,\xi)>-\pi$ and $\bar{q}(t,\xi) = 0$ if 
$v(t,\xi)\leqslant-\pi$.
The system (\ref{ODE3}) is the key tool to prove the following theorem.
\begin{thm}\label{existence2}
Let $u_0 \in \dot{H}^1(\mathds{R}) \cap L^\infty(\mathds{R})$. If 
there exist a Lipschitz function $\phi$ such that $\phi' \in 
L^1(\mathds{R})$ and with $u_0-\phi \in H^1(\mathds{R})$, then there exists 
a global dissipative solution $u$  of the equation \eqref{RB_2}, such 
that $u(t, \cdot)-\phi \in H^1(\mathds{R}) $ for all $t>0$. In addition, for all $t>0$
\begin{equation}\label{Oleinik}
u_x(t,x)\ \leqslant\ 2\,/\,t \qquad (t,x) \in \mathds{R}^+ \times \mathds{R},
\end{equation}
and if $u_0 \in H^1$, then for almost all $t>0$
\begin{equation}\label{lossenergy}
\int_\mathds{R}\left[\,u(t,x)^2\, +\, \ell^2\, u_x(t,x)^2\,\right] \mathrm{d}x\ \leqslant\ 
\int_\mathds{R} \left[\,u_0(x)^2\, +\, \ell^2\, u_0'(x)^2\,\right] \mathrm{d}x.
\end{equation}
\end{thm}

\begin{rem} 
Due to the loss of the Oleinik inequality (cf. Remark \ref{remAntiOl}), 
the system \eqref{ODE} is slightly modified to \eqref{ODE3} in order to obtain {\em dissipative\/} 
solutions of $\mathrm{rB}$ that satisfies the Oleinik inequality 
\eqref{Oleinik}.
\end{rem}

\begin{rem}
In general, if the initial datum satisfies $u_0' \leqslant M \in 
\mathds{R} \cup \{ +\infty \} $, then the Oleinik inequality 
\eqref{Oleinik} can be improved as
\begin{equation}\label{Olgen}
u_x(t,x)\ \leqslant\ 2\,M\,/\,(Mt+2) \qquad (t,x) \in 
\mathds{R}^+ \times \mathds{R},
\end{equation}
as shown in \eqref{Ol_proof} below.
\end{rem}

\textit{Proof of Theorem \ref{existence2}.}
The idea of the proof is similar to Theorem \ref{existence} above and it  
is done in the following steps: 
\begin{itemize}
\item In the first step we prove the existence of the global solution as in Theorem \ref{existence}.
\item In the second step we prove the dissipation of the energy and the Oleinik inequality.
\end{itemize}

{\bf Step 1: Existence of a solution.}
As in the proof of Theorem \ref{existence}, it suffices to 
show that \eqref{ODE3b} and \eqref{ODE3c} are locally well posed in the 
domain $\mathfrak{D} \subset X$, $\mathfrak{D}$ being defined below and 
$X\eqdef\mathcal{C}([0,T],L^\infty(\mathds{R},\mathds{R}^2))$.

Note that if $v$ is near $-\pi$ the right-hand side of \eqref{ODE3b} is 
discontinuous. To avoid this discontinuity, the system \eqref{ODE3} is 
replaced, as in \cite{BressanConstantin2007b}, by

\begin{equation}\label{ODE4}
U_t(t,\xi)\ =\ F(U(t,\xi))\ +\ G(\xi,U(t,\cdot)), \qquad U\ =\ (v,q),
\end{equation}
with
\begin{align*}
F(U )\ \eqdef\  \begin{cases}
      (-\sin^2 \frac{v}{2}, \half\, q\, \sin v)   & v>-\pi, \\
       (-1,0)  & v \leqslant -\pi,
    \end{cases}    \quad
G(U )\ \eqdef \begin{cases}
      (-\, P\, (1\ +\ \cos v), -\, P\, q\, \sin v)   & v>-\pi, \\
       (0,0)  & v \leqslant -\pi ,
    \end{cases}
\end{align*}

Note also that, as long as the solution to \eqref{ODE4} is well defined, replacing 
$v$ by $\max \{ -\pi , v \}$ gives a solution of the equations 
\eqref{ODE3b} and \eqref{ODE3c}.
In the rest of this step, our aim is to show that the system 
\eqref{ODE4} is locally well-posed.
Let $\delta \in ]0, {\textstyle \frac{2 \pi}{3}}]$ and let $\Lambda$ be defined by
\begin{equation}
\Lambda\ \eqdef\ \big\{\, \xi,\ v_0(\xi)\, \in\, ]-\pi,\ \delta-\pi]\, \big\}.
\end{equation}
The equation \eqref{ODE4} implies that, if $v\, \in\, ]-\pi,\ \delta-\pi]\ \subset\ 
]-\pi,\ -\frac{\pi}{3}]$, then $v_t\leqslant-\half$.
Let $\mathfrak{D} \subset X$ satisfy $U(0,\xi)=U_0(\xi)$ and
\begin{subequations}
\begin{align}
1\, /\, C\ \leqslant\ q(t,\xi)\ &\leqslant\ C\ \quad &\forall (t,\xi) \in 
[0,T] \times \mathds{R},\\
\left| \left\{\, \xi,\ \sin^2  \left( v(t,\xi)\, /\, 2 \right) \ \geqslant\ 
{\textstyle \half}  \right\} \right|\ & \leqslant\ C\ \quad &\forall t \in 
[0,T],\\ \label{ULip}
\left\| \, U(t)\ -\ U(s)\, \right\|_\infty\  &\leqslant\ C\, |t-s| \quad &\forall t,s \in [0,T], \\ 
\label{vdif}
v(t,\xi)\ -\ v(s,\xi)\  &\leqslant\ -\ {\textstyle \frac{t-s}{2}}\ \quad 
&\forall  \xi \in \Lambda,\ 0 \leqslant s \leqslant t \leqslant T.
\end{align}
\end{subequations}
Taking $(v,q) \in \mathfrak{D}$ and using \eqref{P,P_x2}, one gets that 
the right-hand sides of \eqref{ODE3b} and \eqref{ODE3c} are bounded. 
However, the inequality \eqref{PLip} is no longer true and we have instead
\begin{equation}
\left\|\, P(U)\ -\ P(\tilde{U}) \right\|_\infty\
\lesssim\  \left\|\, U\ -\ \tilde{U} \right\|_\infty\ +\ \left| \left\{\, \xi,\ 
(v(\xi) + \pi)(\tilde{v}(\xi) + \pi)\ <\ 0 \right\} \right| ,
\end{equation}
which implies that
\begin{equation}
\left\|\, F(U)\ -\ F(\tilde{U}))\right\|_\infty\ \lesssim\  \left\|U\ -\ \tilde{U} \right\|_\infty ,
\end{equation}
\begin{equation}
\left\|\, G(U)\ -\ G(\tilde{U})) \right\|_\infty\ \lesssim\ \left\|U\ -\ \tilde{U} \right\|_\infty\ +\ \Big| 
\Big\{ \xi,\ (v(\xi)\ +\ \pi)(\tilde{v}(\xi)\ +\ \pi)\ <\ 0 \Big\} \Big|.
\end{equation}
In order to estimate the second term of the right-hand side of the last equation,  
the crossing time $\tau$ is defined as
\begin{equation}\label{crossingtime}
\tau(\xi)\ \eqdef\ \sup\ \{ t\, \in\, [0,T],\ v(t,\xi)\ >\ -\, \pi \}.
\end{equation}
Note that the equation \eqref{ULip} implies that $|v(t,\xi) - v_0(\xi)| 
\leqslant C\, t$. So, if $\xi\, \notin\, \Lambda$ then
\begin{equation*}
\min\left\{ \tau(\xi),\  \tilde{\tau}(\xi) \right\}\ \geqslant\ 
\delta\, /\, C.
\end{equation*}
Taking $T$ small enough ($T < \delta/C$) and using 
\eqref{vdif}, one obtains
\begin{align*}
\int_0^T \Big| \Big\{ \xi,\ (v(\tau,\xi)\ +\ \pi)(\tilde{v}(\tau,\xi)\ +\ \pi)\ <\ 0 
\Big\} \Big|\ \mathrm{d} \tau\
& \leqslant\  \int_\Lambda |\tau(\xi)\ -\ \tilde{\tau}(\xi)|\ \mathrm{d} \xi 
\\
& \leqslant 2\, |\, \Lambda\, |\, \left\|\, U\ -\ \tilde{U} \right\|_\infty .
\end{align*}
Now, the Picard operator
\begin{equation}
(\mathcal{P}(U))(t,\xi)\ =\ U_0\ +\ \int_0^t \left[ F(U)\ +\ G(U) \right]\, 
\mathrm{d} \tau,
\end{equation}
satisfies
\begin{equation}
\left\|\, \mathcal{P}(U)\ -\ \mathcal{P}(\tilde{U}) \right\|_\infty\ \leqslant\ K\,
(T\ +\ |\, \Lambda\, |) \left\|U\ -\ \tilde{U} \right\|_\infty ,
\end{equation}
where $K$ depends only on $C$ and $\ell$.
Since $\sin^2 {\textstyle \frac{v_0}{2}} \in L^1$, by choosing $\delta>0$ small enough, 
one can make $|\Lambda|$ arbitrary small. Choosing also $T$ small enough, one obtains 
the local existence of the solution of the system \eqref{ODE4}, yielding a solution of 
\eqref{ODE3}. The rest of the proof of the existence can be done following the proof of 
Theorem \ref{existence}.\\

{\bf Step 2: Oleinik inequality and the dissipation of the energy.}
The equation \eqref{ODE3b} implies that if $v(0,\xi)\leqslant 0$, then 
for all $t \geqslant 0$ $v(t,\xi)$ remains in $[-\pi,0]$. If  $v_0(\xi) 
\in ]0,\pi[$ then, as long as $v$ is positive, the following inequality 
holds
\begin{align*}
\left[ \arctan \frac{v}{2} \right]_t \leqslant - \half \arctan^2 
\frac{v}{2}.
\end{align*}
This implies that, if $\arctan \frac{v_0(\xi)}{2} \leqslant M$, then
\begin{equation}\label{Ol_proof}
  u_x=\arctan \frac{v(t,\xi)}{2} \leqslant \frac{2M}{Mt+2}.
\end{equation}
The Oleinik inequality \eqref{Oleinik} follows taking $M=+\infty$ and using 
\eqref{u_x}.

In order to prove the dissipation of the energy \eqref{localdiss}, let $\psi$ 
be a non-negative test function, then we follow the same computations in the proof of Theorem 
\ref{existence}. Since \eqref{measure} is no longer true for the system \eqref{ODE3}, 
one can obtain from \eqref{localene} that 
\begin{gather*}
\iint\limits_{[0,+\infty[ \times \mathds{R}}\! \left[- u_x^2 \psi_t\, -\, u u_x^2 \psi_x\, +\, u_x P \right]\, \mathrm{d}t\, \mathrm{d}x\, -\! \int\limits_\mathds{R}\! u_0'^2(x) \psi(0,x)\, \mathrm{d}x\\
= -\! \int\limits_{\{\tau(\xi)< +\infty\}}\! q(\tau(\xi), \xi) \tilde{\psi}(\tau(\xi), \xi)\, 
\mathrm{d}\xi\ \leqslant\ 0,
\end{gather*}
where $\tau(\xi)$ is the crossing time defined as $\tau(\xi) \eqdef \sup\ \{ t\/ 
\geqslant\/ 0,\/ v(t,\xi) > -\/\pi \}$. Since \eqref{rBene1} is satisfied (see Remark \ref{energyeq}), the dissipation of the energy \eqref{localdiss} follows.

If $u_0 \in H^1$, as in the last step of the proof of Theorem \ref{existence}, one can 
show that \eqref{conservation2} and \eqref{energyrelation} hold for the 
solution of \eqref{ODE3}, while the measure in \eqref{measure} is not 
always zero. Then, the dissipation of the energy \eqref{lossenergy} follows.
\qed

\section{Traveling waves of permanent form}\label{s:waves}

The rB equation \eqref{rB0} can be written in the conservative form
\begin{equation}\label{RBmom}
\left[\,u\,-\,\ell^2\,u_{xx}\,\right]_t\ +\,\left[\, \half\, u^2\,-\,\ell^2\,u\,u_{xx}\,
-\,\half\,\ell^2\, u_x^{\,2}\,\right]_x\ =\ 0.
\end{equation}
In this section, we describe traveling waves of permanent form, i.e., we seek bounded weak solutions 
of \eqref{RBmom} having the form $u=u(x-ct)$. 
We find a great variety of such waves, roughly comprising a subset of the plethora of traveling wave solutions 
of the Camassa-Holm equation found by Lenells~\cite{Lenells2005}.  For CH, Lenells found 
peakons, cuspons, stumpons, smooth periodic waves, monotone waves, and composite waves of various kinds.
When we consider the requirements imposed by energy conservation or dissipation, however, 
many composite constructions are eliminated.
In particular, the only
traveling waves of rB we find that are dissipative in the sense of Definition~\ref{defdis}
are monotone weakly singular shock layers that correspond to
entropy-satisfying shocks of the inviscid Burgers equation in the limit $\ell\to0$.
The rB equation has the nice property that such shock layers exist for every 
entropy-satisfying shock. 
This is the scalar analog of the property found by \citet{PuEtAl2018} for the  
nondispersively regularized shallow water system~\eqref{rSV1}--\eqref{rSV3}.

Since the $\mathrm{rB}$ equation is Galilean invariant, we can work 
in the frame of reference moving with the wave, where the motion is steady. 
Moreover, by the rescaling $x/\ell \to x$ we can presume $\ell=1$.
I.e., we look for (weak) solutions such that $u=u(x/\ell)$.

For steady motions with $\ell=1$, the momentum flux $S$, given by 
\begin{equation}
S\ \eqdef\ \half\,u^2 - u\,u_{xx} - \half\,u_x^{\,2} \ = \ \half\,u^2 + P,\label{defS}
\end{equation}
is constant, cf.~\eqref{Ric}.  In any open set where $u\ne0$, equation \eqref{defS} is an ODE 
and any weak solution must be smooth.
The energy conservation law \eqref{rBene0} then implies that the energy flux $F$, given by 
\begin{align}
F(x)\ &\eqdef\ \third\,u^3- u^2\,u_{xx} \ =\  \third\,u^3 +u\,P +\half \,u \,u_x^2, \label{defF}
\end{align}
is locally constant, cf.~\eqref{rBene}.
Eliminating $u_{xx}$ between (\ref{defS}) and (\ref{defF}), one obtains the first-order 
differential equation
\begin{equation}\label{TWSF}
\half\,u\,u_x^{\,2}\ =\ F\ -\ S\,u\ +\ \sixth u^3.
\end{equation}

\subsection{Local analysis of weak solutions} 
As it turns out, the energy flux $F$ may be discontinuous at points where $u=0$.
Consider the possibility of a singularity at $x=x_0$, which can exist only if $u(x_0)=0$. Equation \eqref{TWSF} then yields 
the following asymptotic behavior: As $x\to x_0$ (from either the right or the left),
\begin{align}
 u\, u_x^2\ \sim\ 2\,F 
 &\quad\mbox{ if  }Fu\,>\,0, \\  
 u_x^2\ \sim\ -2\,S 
 &\quad\mbox{ if  }F\, =\, 0 \text{ and } S\, <\, 0,
\end{align}
implying that, as $x\to x_0$,
\begin{align}
\label{a:cusp}
|u|\ \sim\  \left(\threehalf\, \sqrt{|2F|}\, |x-x_0|\right)^{2/3}
\qquad &\text{if } Fu\,>\,0, \\
\label{a:peak}
|u|\ \sim\ \sqrt{|2S|}\ |x-x_0| \quad\qquad\qquad &\text{if }F\,=\,0\text{ and } S\,<\,0.
\end{align}
By consequence, the quantity $u_x^2 = O(|x-x_0|^{-2/3})$ is locally integrable.
Hence by \eqref{defS}, one can construct a valid weak solution of \eqref{RBmom}
by solving \eqref{TWSF} separately for $x<x_0$ and $x>x_0$,
allowing $F$ to jump between any two real values from left to right across $x_0$,

\subsection{Cuspons and periodic cuspons}
In the simplest case, $F$ is globally a nonzero constant.
This is necessary for the wave to be a conservative solution according to Definition~\ref{defcon},
since \eqref{rBene} requires $F_x=0$.  Then we obtain a family of waves with cusps
having the behavior in \eqref{a:cusp} for $x_0$ either at a single point or at any point in a periodic grid.

Note that $u(-x)$ is a solution of \eqref{TWSF} whenever $u(x)$ is, and $-u$ is a solution for $-F$ in place of $F$.
Without loss, then, we may consider the case when $u>0$ for $x>0=x_0$, and $F>0$.
Necessarily, if $u$ is to be bounded, the cubic polynomial on the right-hand side of \eqref{TWSF} must have two positive roots
at points $u_1\geqslant u_0>0$ and a negative root at $-u_1-u_0$, related to $F$ and $S$ via
\begin{equation}
F = \sixth u_0u_1(u_0+u_1), \qquad S = \sixth (u_0^2+u_0u_1+u_1^2).
\end{equation}
By consequence, we can separate variables in \eqref{TWSF} and find that $u = \eta(x)$ where
$\eta$ is determined implicitly by $x$ from the relation
\begin{equation}
H(\eta) \eqdef \int_0^\eta \left(\frac{3v}{(u_0-v)(u_1-v)(u_0+u_1+v)}\right)^{1/2}\,dv = x, \qquad 0<x<x_*=H(u_0)\le\infty.
\end{equation}

We obtain a {\em periodic cuspon} solution whenever $u_0<u_1$. In this case the integral converges at $\eta=u_0$ to a finite value,
and $u(x)$ can be defined by reflection 
about $x_0=0$ as 
\begin{equation}
u(x) = \eta(|x|/\ell) ,\qquad |x|\leqslant x_*\ell, 
\end{equation}
then extended as smooth and periodic with period $2x_*\ell$ and maximum value $u_0=u(x_*\ell)$. 
Here we have put back the scale parameter $\ell$ to obtain a stationary weak solution of \eqref{RBmom} 
valid for any $\ell>0$. 

We find a {\em cuspon} with $u(x) \to u_0$ as $x\to \pm\infty$, 
provided by the same formulas in the case $u_0=u_1$, when we find $x_*=H(u_0)=+\infty$.
This has a single singular point at $x_0=0$.
In this case, 
\begin{equation}\label{e:cuspon}
F = \third u_0^3, \qquad S = \half u_0^2, \qquad 
H(\eta) = \int_0^\eta \left(\frac{3v}{2u_0+v}\right)^{1/2}\frac{dv}{u_0-v} =|x|.
\end{equation}

A similar analysis provides negative cuspons $u<0$ with reversed signs $F<0$ and $u_1\leqslant u_0<0$. 
Necessarily $S>0$ in this case also.
In case $F=0$ globally, one has only the trivial solution $u\equiv0$, 
for there are no other bounded solutions, since $u_{xx}=\third u$ for $u\ne0$ by \eqref{defF}.
Finally, also we note that 

In summary we can state the following, noting that $\min_{u>0}(|F|-Su+\sixth u^3)\le0$ in all cases.

\begin{pro} For globally constant energy flux $F$,
a bounded nonzero stationary solution $u$ of \eqref{RBmom} exists, having the form 
\begin{equation}
u(x) = \eta(|x|/\ell) \sgn F ,\qquad \mbox{for }|x|\leqslant x_*\ell, 
\end{equation}
provided
\begin{equation}\label{c:SF}
0< 3|F| \leqslant (2S)^{3/2},
\end{equation}
where we have equality for cupsons, strict inequality for periodic cuspons.  
All these cuspons and periodic cuspons are conservative solutions according to Definition~\ref{defcon}.
\end{pro}

\subsection{Shock layers}
For a stationary solution $u$ to be {\em dissipative} according to Definition~\ref{defdis}, we require $u_x\le0$
to satisfy the Oleinik inequality~\eqref{oleinik}. 
The energy flux $F(x)$ can be discontinuous, but it must be {\em decreasing} across singularities in order to satisfy 
the energy dissipation inequality~\eqref{localdiss}.

Thus we require $u>0$ and $F=F_->0$ for $x<0=x_0$, while $u<0$ and $F=F_+<0$ for $x>0$. 
Since $u(x)$ is required to be bounded and monotone we should have
\begin{equation}
u(x) \to \begin{cases}
u_- & \mbox{as \ $x\to-\infty$},\\
u_+ & \mbox{as \ $x\to+\infty$},
\end{cases}
\end{equation}  
where necessarily $S>0$ and $u_-=-u_+ = \sqrt{2S}$ by taking $x\to\pm\infty$ in \eqref{defS}.
Since we must solve the same ODE as in the previous subsection, to obtain a global monotonic solution
it is necessary that 
\begin{equation}
|u_\pm|= (2S)^{1/2} = (3|F_\pm|)^{1/3}.
\end{equation}
We obtain such solutions by taking $u_0=u_-$ in the formulas in \eqref{e:cuspon},
and taking the odd extension of the left half of the cuspon, setting
\begin{equation}
 u(x) = \begin{cases} \eta(-x/\ell), & x<0,\\ -\eta(x/\ell), & x>0.\end{cases}
\end{equation}

As before this yields a solution of \eqref{RBmom} for any $\ell>0$. 
In the limit $\ell\to0$ we evidently obtain 
any stationary entropy-satifying shock for the inviscid Burgers equation, 
which must take two values $u_->u_+$ with $u_-^2=u_+^2$.
By Galilean boosts we obtain traveling weakly singular shock layer solutions which converge in the limit $\ell\to0$
to any arbitrary entropy-satisfying simple shock for the inviscid Burgers equation.
Thereby we obtain the following.

\begin{pro}
Corresponding to any entropy-satisfying simple shock for the inviscid Burgers equation, 
taking constant values $u_->u_+$ respectively for $x<ct$ and $x>ct$ with $c=\half(u_-+u_+)$, 
the regularized Burgers equation \eqref{RBmom} admits a weakly singular shock layer solution satsifying
\begin{equation}
u(x,t) = c -\eta(|x-ct|/\ell)\sgn(x-ct)
\to \begin{cases} u_- & x\to-\infty,\\ u_+ & x\to+\infty,
\end{cases}
\end{equation}
where $\eta(x)$ is determined from \eqref{e:cuspon} with $u_0=\half(u_--u_+)$.
\end{pro}

\begin{rem}
We note from \eqref{localdiss} that the rate of energy dissipation for the stationary shock layer is 
\[
F_- - F_+ = \twothird u_-^3  = {\textstyle\frac1{12}} (u_- - u_+)^3 , 
\]
and all of the energy dissipation occurs at the location of the weak singularity.
This rate is exactly the same as for the corresponding inviscid Burgers shock, 
which is famously cubic in the amplitude of the shock. 
\end{rem}

\begin{figure}
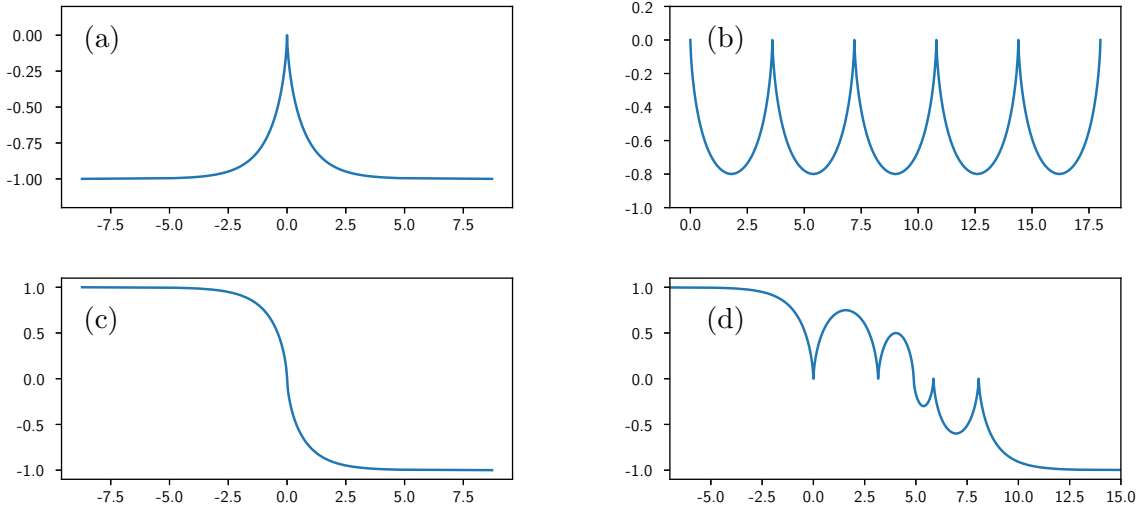

\begin{tabular}{cc}
\scalebox{0.62}{\input{fig1a-cuspon.pgf}} &
\scalebox{0.62}{\input{fig1b-periodic.pgf}} \\
\scalebox{0.62}{\input{fig1c-shock.pgf}} &
\scalebox{0.62}{\input{fig1d-composite.pgf}}
\end{tabular}
\put(-415,78){(a)}
\put(-182,78){(b)}
\put(-415,-28){(c)}
\put(-182,-28){(d)}
\caption{Types of weakly singular stationary waves:
(a) cuspon; (b) periodic cuspon; (c) shock layer; (d) composite wave}
\end{figure}

\subsection{Composite waves}  If one discards the criteria we have imposed to
find conservative or dissipative solutions, a great many more stationary weak
solutions of \eqref{RBmom} can be constructed, by joining together segments of (periodic) cuspons between consecutive zeros, 
while allowing the energy flux $F(x)$ to jump discontinuously at the zeros, 
in an essentially arbitrary way subject only to the inequality \eqref{c:SF}. 

We will not describe these solutions in any further detail, since we have already described
the only such solutions that are conservative or dissipative according to Definitions~\ref{defcon} and \ref{defdis}.
At any isolated singular point $x_0$ where the left-to-right jump $F(x_0+)-F(x_0-)$ is negative, energy is {\em dissipated},
while energy is {\em generated} if this jump is positive.  
One can construct composite solutions having multiple zeros that all dissipate energy, 
by taking $F(x)$ piecewise constant and non-increasing.
One such composite solution is plotted in Figure~1(d).
These are not dissipative solutions in the sense of Definition~\ref{defdis},
however, since they are not monotonic hence violate the Oleinik inequality.

\subsection{Non-existing waves}
Some of the kinds of waves found by Lenells for CH are not possible for rB.  
In particular, we find that rB does not admit smooth periodic waves, 
peakons with bounded derivatives, or ``stumpons'': 
Such waves can be made stationary by a Galilean transformation, 
then must have a single sign, which may be assumed nonnegative.  
If $\min u=0$ we must have $S>0$ and the wave must be a cuspon, as found above. 
And if $\min u>0$, then the cubic polynomial in \eqref{TWSF} must
take positive values between two zeros $u_0=\min u<\max u=u_1$. 
But this is not possible since $-u_0-u_1$ is the only other zero.

\section{The limiting cases $\ell\to0$ and $\ell\to+\infty$ for dissipative solutions }\label{Compactness}

Taking formally $\ell = 0$, the $\mathrm{rB}$ equation becomes the classical Burgers equation, 
and letting $\ell\to+\infty$ it becomes the Hunter--Saxton equation. In this section, we study 
the compactness of the dissipative solutions when taking $\ell\to 0$ and $\ell \to +\infty$. 

Let the initial datum $u_0$ be taken in $H^1$, with $u_0'\in L^1$ and
$M\eqdef\sup_{x \in \mathds{R}} u_0'(x) < +\infty$. Let also $u^\ell$ be the 
dissipative solution of the $\mathrm{rB}$ equation given in Theorem 
\ref{existence2}. 
In order to take the limit, an estimate on the total variation of $u^\ell$ is needed. 
For that purpose, the following Lemma is given
\begin{lem}{\textbf{[BV estimate]}}\label{TV}
If $u_0'$ satisfies the conditions of Theorem \ref{existence2} with $u_0' \in L^1$ and $u_0'(x) \leqslant M\ 
\forall x$, then for all $t \in \mathds{R}^+$
\begin{equation}
\mathrm{TV} u^\ell(t,\cdot)\ =\ \left\| u^\ell_x(t,\cdot) \right\|_1\ \leqslant\ \left\|u_0'\right\|_1\, \left( \frac{Mt+2}{2} \right)^2.
\end{equation}
\end{lem}

\proof 
For $v \in 
]-\pi,\pi[$, the equation \eqref{u_xi} implies
\begin{equation}
\tilde{s} \eqdef 
\mathrm{sgn}(u^\ell_\xi)\ =\ \mathrm{sgn}(\sin(v^\ell))\ =\ \mathrm{sgn}(\sin(\frac{v^\ell}{2})), 
\qquad \cos (\frac{v^\ell}{2})\ \geqslant\ 0.
\end{equation}
Note that $\tan (v^\ell/2) \leqslant \tan (v_0/2) =  2M/(Mt+2)$ from \eqref{Ol_proof}. 
Differentiating \eqref{ODE3a} w.r.t $\xi$, multiplying by $\tilde{s}$ --- and using \eqref{P_xi}, 
\eqref{u_xi} and $\sin v = 2 \sin (v/2) \cos (v/2) $ --- one gets
\begin{align}
\frac{\mathrm{d}}{\mathrm{d} t} \int_\mathds{R} |u^\ell_\xi|\ \mathrm{d} \xi\ &=\ -\ell^2 
\int_\mathds{R} \tilde{s}\ (P_x)_\xi\  \mathrm{d} \xi
=\ -\ell^2 \int_{\{\tilde{s}\ >\ 0\}} (P_x)_\xi\   \mathrm{d} \xi\ +\ \ell^2\,  
\int_{\{\tilde{s}\ <\ 0\}} (P_x)_\xi\   \mathrm{d} \xi \nonumber \\
&=\ -2\, \ell^2\,  \int_{\{\tilde{s}\ >\ 0\}} (P_x)_\xi\   \mathrm{d} \xi\ +\  
\ell^2\, \int_\mathds{R} (P_x)_\xi\ \mathrm{d} \xi \nonumber \\
&=\ -2\  \int_{\{\tilde{s}\ >\
 0\}} \left( q^\ell\, P\, \cos^2 
\frac{v^\ell}{2}\ -\ \half q^\ell\, \sin^2 \frac{v^\ell}{2} \right)\,   
\mathrm{d} \xi \nonumber \\
& \leqslant\  \int_{\{\tilde{s}\ >\ 0\}} q^\ell\, \sin \frac{v^\ell}{2}\, \cos 
\frac{v^\ell}{2}\, \tan \frac{v^\ell}{2}\ \mathrm{d} \xi \nonumber \\
& \leqslant\  \frac{2M }{Mt+2}\, \int_\mathds{R} |u^\ell_\xi|\ \mathrm{d} \xi .
\end{align}
Gronwall lemma then implies that
\begin{equation}
\|u_\xi\|_1\ \leqslant\ \|(u_0)_\xi\|_1\, \left( \frac{Mt+2}{2} \right)^2.
\end{equation}
Note that the last inequality is on the $\xi-$variable. Using that the application 
$\xi \mapsto y(t,\xi)$ is not decreasing for all $t$ and using that 
$\mathrm{TV} f = \|f'\|_1$ for smooth solutions ($ f\in W^{1,1}_{loc}$),  
the result follows. \qed 

\subsection{The limiting case $\ell \to 0$}

The goal of this subsection is to show that when $\ell\to0$, the dissipative solution 
$u^\ell$ converges {(up to a subsequence)} to a function $u$ satisfying the Burgers equation 
with a source term: 
\begin{equation}\label{limRB}
u_t\, +\, \half\left[u^2\right]_x\, =\, -\, \mu_x,
\end{equation}
where $\mu$ is a measure such that $0 \leqslant \mu \in L^\infty ([0,+\infty[, \mathcal{M}^1)$. 
In Proposition \ref{mu} below, we show that the measure $\mu$ is zero before the appearance 
of singularities. The question whether or not $\mu$ is zero after singularities is open, in general.
The following theorem can be stated
\begin{thm}\label{ThmCompactness}
Let $u_0 \in H^1$, such that $u_0' \in L^1$ and $u_0'(x) \leqslant M\ 
\forall x$, then there exists $u \in  L^\infty([0,T], 
BV(\mathds{R}))$ for all $T>0$, such 
that there exists a subsequence of $u^\ell$ (also noted $u^\ell$) and for all interval $\mathcal{I} \subset \mathds{R}$ we have 
\begin{equation}\label{convergence0}
u^\ell  \xrightarrow{\ell \to 0}\,  u\ \text{in}\ 
\mathcal{C}([0,T],L^1(\mathcal{I})), 
\end{equation}
and $u$ satisfies the equation \eqref{limRB}.
Moreover, $u$ satisfies the Oleinik inequality
\begin{equation}
u_x(t,x) \leqslant \frac{2M     }{Mt+2} \quad \textrm{ in } 
\mathcal{D}'(\mathds{R}).
\end{equation}
\end{thm}
\begin{rem}
If $\mu=0$ then, due to the Oleinik inequality, $u$ is the unique entropy solution of 
the Burgers equation.
\end{rem}
In order to prove Theorem \ref{ThmCompactness}, the following definition and lemma are needed: 

\noindent
Let $\mathcal{I} \subset \mathds{R}$ be a bounded interval and let 
\begin{equation}\label{Wdef}
W(\mathcal{I})\ \eqdef\ \left\{ f\, \in\, \mathcal{D}'(\mathcal{I}),\ \exists\,
F\, \in\, L^1(\mathcal{I})\ \textrm{such that } F'\ =\ f \right\} ,
\end{equation}
where the norm of the space $W(\mathcal{I})$ is given by
\begin{equation}\label{norm}
\|f\|_{W(\mathcal{I})}\ \eqdef\ \inf_{c\, \in\, \mathds{R}} 
\|F + c\,\|_{L^1(\mathcal{I})}\ =\ \min_{c\, \in\, \mathds{R}} 
\|F + c\,\|_{L^1(\mathcal{I})} .
\end{equation}
\begin{lem} The space $W(\mathcal{I})$ is a Banach space and the 
embedding
\begin{equation}
L^1(\mathcal{I})\ \hookrightarrow\ W(\mathcal{I}),
\end{equation}
is continuous.
\end{lem}
\proof Let $(f_n)_{n\in\mathds{N}}$ be a Cauchy sequence in $W(\mathcal{I})$ and let 
$F_n$ be a primitive of $f_n$. From the definition of the norm \eqref{norm}, there exists 
a constant $c_n$ such that $(\tilde{F}_n-c_n)_{n\in\mathds{N}}$ (where $\tilde{F}_n=F_n
+c_n$) is a Cauchy sequence in $L^1(\mathcal{I})$. 
Let $\tilde{F}$ be the limit of $\tilde{F}_n$ in $L^1(\mathcal{I})$. Then
\begin{equation}
\|f_n\ -\ \tilde{F}'\|_{W(\mathcal{I})}\ \leqslant\ 
\|\tilde{F_n}\ -\ \tilde{F} \|_{L^1(\mathcal{I})},
\end{equation}
implying that $W(\mathcal{I})$ is a Banach space.
Now, the continuous embedding can be proved.

If $f \in L^1(\mathcal{I})$, then $F(x) - F(a) = \int_a^x f(y)\, \mathrm{d}y$  
for almost all $x,a \in \mathcal{I}$. Therefore, 
\begin{equation}
\|f\|_{W(\mathcal{I})}\ \leqslant\ \int_\mathcal{I} |F(x)\ -\ F(a)|\ 
\mathrm{d}x\ \leqslant\ |\mathcal{I}| \int_\mathcal{I} |f(y)|\ \mathrm{d}y 
,
\end{equation}
which ends the proof. \qed

The previous lemma and Helly's selection theorem imply that
\begin{equation}
W^{1,1}(\mathcal{I})\ \hookrightarrow\ L^1(\mathcal{I})\ \hookrightarrow\ 
W(\mathcal{I}),
\end{equation}
where the first embedding is compact and the second is continuous. \\

\noindent
\textbf{Proof of Theorem \ref{ThmCompactness}:}
Let the compact set $[0,T] \times \mathcal{I} \subset \mathds{R}^+ \times 
\mathds{R}$. Supposing that $\ell \leqslant 1$ then, from \eqref{lossenergy}, 
the dissipative solutions of $\mathrm{rB}$ satisfies
\begin{equation}\label{H^1}
\|u^\ell\|_2^2\ \leqslant\ \|u_0\|_{H^1}^2, \qquad \qquad \ell^2 \|P\|_1\ =\ 
\half\, \ell^2\, \|u^\ell_x\|_2^2\ \leqslant\ \half\, \|u_0\|_{H^1}^2,
\end{equation}
implying that $u^\ell$ is uniformly bounded on $L^\infty ([0,T], 
L^2(\mathds{R}))$. Subsequently, it is also uniformly bounded on $L^\infty ([0,T], 
L^1(\mathcal{I}))$. Because Lemma \ref{TV} yields that $u^\ell$ is bounded on $L^\infty([0,T],
W^{1,1}(\mathcal{I}))$, and the equation \eqref{H^1} implies that $\half {u^\ell}^{\,2} 
+ \ell^2 P$ is uniformly bounded on $L^\infty([0,T],L^1(\mathcal{I}))$, then since 
$u^\ell_t = -\left( \half {u^\ell}^{\,2} + \ell^2 P \right)_x$, \eqref{norm} implies 
that $u^\ell_t$ is bounded on $L^\infty([0,T],W(\mathcal{I}))$. 
Then, using the Aubin theorem 
(see Corollary 4 in \cite{Simon}), the compactness follows.

The quantity $\ell^2\/P$ is non-negative and bounded in $L^\infty([0,+\infty[,L^1(\mathds{R}))$, 
implying the existence of a non-negative measure $\mu \in L^\infty([0,+\infty[,\mathcal{M}^1(
\mathds{R}))$ such that $\ell^2\/P$ converges {(up to a subsequence)} weakly to $\mu$.
The equation \eqref{limRB} follows taking the limit $\ell \to 0$ in the weak formulation of 
\eqref{RB_2}. Finally, taking the limit in the weak formulation of \eqref{Olgen}, we 
can prove that $u_x(t,x) \leqslant 2M/(Mt+2)$. \qed

The question whether or not $\mu = 0$ is open. The following proposition shows that when 
$\ell\to 0$ for smooth solutions (i.e., before appearance of singularities), $u^\ell$ converges 
to the unique solution $u$ of the classical Burgers equation.
\begin{pro}\label{mu}
If $u_0$ is in $H^s \cap BV$ with $s \geqslant 3$, then for $t < {1/ \sup_x |u_0'(x)|}$ 
we have
\begin{equation}
\mu\ =\ 0.
\end{equation}
\end{pro}
\proof
From Theorem \ref{InteractionTime} {and Remark \ref{BBT}}, we can find a uniform upper 
bound on $u^\ell_x$ in the space $L^\infty([0,T],L^\infty(\mathds{R}))$ with 
$T<1/ \sup_x |u_0'(x)|$, which implies that $\ell^2\/P \to 0$. \qed

\subsection{The limiting case $\ell \to +\infty$}

The goal of this subsection is to show that, when $\ell 
\to +\infty$, the dissipative solution $u^\ell$ converges {(up to a subsequence)} to a function 
$u$ that satisfies:
\begin{equation}\label{limHS}
\left[ u_t\, +\, \half(u^2)_x \right]_x \, =\, \nu,
\end{equation}
where $0 \leqslant \nu \in L^\infty ([0,+\infty[, \mathcal{M}^1)$. In Proposition 
\ref{nu} below, we show that before the appearance of singularities, the measure $\nu =  
u_x^2/2$. The question whether or not $\nu =  u_x^2/2$ in general is posed.
We have the following theorem: 
\begin{thm}\label{ThmCompactness2}
Let $u_0 \in H^1$ such that $u_0' \in L^1$ and $u_0'(x) \leqslant M\ 
\forall x$, then there exists $u \in 
 L^\infty([0,T], 
BV(\mathds{R}))$ for all $T>0$, such 
that there exists a subsequence of $u^\ell$ (noted also $u^\ell$) and for all interval $\mathcal{I}  \subset \mathds{R}$ we have
\begin{equation}\label{convergence}
u^\ell \xrightarrow{\ell \to +\infty} u\ \text{in}\ 
\mathcal{C}([0,T],L^1(\mathcal{I})),
\end{equation}
and $u$ satisfies the equation \eqref{limHS}.
Moreover, $u$ satisfies the Oleinik inequality
\begin{equation}
u_x(t,x) \leqslant \frac{2M     }{Mt+2} \quad \textrm{ in } 
\mathcal{D}'(\mathds{R}).
\end{equation}
\end{thm}
\begin{rem}
If $\nu = \half u_x^2$ then $u$ is a dissipative solution of the Hunter--Saxton equation \cite{BressanConstantin05}.
\end{rem}

\proof
Let the compact set $[0,T] \times \mathcal{I} \subset \mathds{R}^+ \times 
\mathds{R}$. Supposing that $\ell \geqslant 1$ then, from \eqref{lossenergy}, 
the dissipative solution of $\mathrm{rB}$ satisfies
\begin{equation}\label{2H^1}
\|P\|_1\ =\ 
\half\,  \|u^\ell_x\|_2^2\ \leqslant\ \half\, \|u_0\|_{H^1}^2.
\end{equation}
Using Lemma \ref{TV}, one gets that 
$u^\ell$ is bounded in $L^\infty ([0,T] \times \mathds{R})$ and 
\begin{equation}
\int_\mathds{R} |u^\ell(t,x+h)\ -\ u^\ell(t,x)|\, \mathrm{d}x\ \leqslant\ \left\|u_0'\right\|_1\, \left( \frac{M\/ T+2}{2} \right)^2 |h|.
\end{equation}
Integrating \eqref{RB_2} between $t_1$ and $t_2$, one obtains
\begin{equation}
u^\ell(t_1,x)\ -\ u^\ell(t_2,x)\ =\ \int_{t_1}^{t_2} \left( u^\ell\, u^\ell_x\ +\ \ell^2\, P_x\, \right)\, \mathrm{d}t.
\end{equation} 
Using Lemma \ref{TV}, inequality \eqref{2H^1} and $$ \|P_x\|_\infty\  \leqslant\ \quat\, \ell^{-2}\, \|u_x^\ell \|_2^2 ,$$ 
we can show that there exists $B\, =\, B(T,\mathcal{I})$ such that 
\begin{equation}
\int_\mathds{\mathcal{I}} |u^\ell(t_2,x)\ -\ u^\ell(t_1,x)|\, \mathrm{d}x\ \leqslant\ B\, |t_2-t_1|.
\end{equation}
The compactness follows using Theorem A.8 in \cite{HoldenRisebro2015}.

The quantity
$\half\, {u^{\ell}_x}^{\,2}$ is non-negative and bounded in $L^\infty([0,+\infty[,L^1(\mathds{R}))$, 
which implies that there exists a non-negative measure $\nu \in L^\infty([0,+\infty[,\mathcal{M}
^1(\mathds{R}))$ such that $P$ converges {(up to a subsequence)} weakly to $\nu$.
The equation \eqref{limHS} follows by taking the limit $\ell \to +\infty$, in 
the weak formulation of \eqref{RB_x}.
Finally, taking the limit in the weak formulation of \eqref{Olgen}, we 
can prove that $u_x(t,x) \leqslant \frac{2M     }{Mt+2}$. \qed

The question whether or not the equality always holds $\nu =  u_x^2/2$ is open. 
The following proposition shows that, when $\ell\to +\infty$ for smooth solutions (before 
appearance of singularities), $u^\ell$ converges to a dissipative solution $u$ of the Hunter--Saxton 
equation \cite{BressanConstantin05}.
\begin{pro}\label{nu}
If $u_0$ is in $H^s \cap BV$ with $s \geqslant 3$, then for $t < {1/ \sup_x |u_0'(x)|}$
we have
\begin{equation}
\nu\ =\ \half\, u_x^2.
\end{equation}
\end{pro}
\proof
From Theorem \ref{InteractionTime} {and Remark \ref{BBT}}, we can find a uniform upper 
bound on $u^\ell_x$ in the space $L^\infty([0,T],L^\infty(\mathds{R}))$ with 
$T<1/ \sup_x |u_0'(x)|$, which implies that the convergence $u_x^\ell$ to $u_x$ is strong. 
Thus, ${u_x^\ell}^{\, 2} \to  u_x^2$. \qed

\section{Optimality of the $\dot{H}_{loc}$   space}\label{secoptimality}

In the previous sections (see Proposition \ref{blowup},  Theorem \ref{InteractionTime} 
and Theorem \ref{existence2}), we have shown, on one side, that even if the initial 
datum $u_0$ is smooth, there exists a finite blow-up time $T^*>0$ such that
\begin{equation}\label{losingreg}
\inf\limits_{x \in \mathds{R}} u_x(t,x)\ >\ -\infty \ \forall t<T^*, \qquad   \inf\limits_{x \in \mathds{R}} u_x(T^*,x)\ =\ -\infty .
\end{equation}
On the other side, the Oleinik inequality \eqref{Oleinik} shows that, even if the initial 
datum is not Lipschitz, the derivative of the solution becomes instantly bounded from above, 
i.e.
\begin{equation}\label{gainingreg}
\sup\limits_{x \in \mathds{R}} u_0'(x)\ =\ +\infty ,\ \qquad   \sup\limits_{x \in \mathds{R}} u_x(t,x)\ <\ +\infty\quad \forall t>0.
\end{equation}
\begin{rem}\label{Riccati_Oleinik}
If the derivative of the initial datum is bounded from below and not from above, it will be 
instantly bounded from both sides \footnote{Note that the gain of regularity \eqref{gainingreg} is instantaneous, while the loss of regularity \eqref{losingreg} needs some time.} and, after $T^*$, 
it will be bounded from above and not from below.
\end{rem}
This remark is important to prove that the space $\dot{H}^1_{loc}$ is the best space to obtain 
global (in time) solutions, the optimality being in the following sense.
\begin{thm}\label{optimality}
Let $\delta>0$ and $g(h) \eqdef \left[ \ln |h| \right]^\delta$, then there exist $u_0 \in H^1 \cap W^{1,\infty}$, $T>0$ and a compact set $\mathcal{K}$, such that there exists a solution $u$ of \eqref{RB_2} satisfying
\begin{equation}
\int_\mathds{R} u_0'(x)^2\, g(u'_0(x))\, \mathrm{d}x\ <\ +\infty, \qquad \qquad \int_\mathcal{K} u_x(T,x)^2\, g(u_x(T,x))\, \mathrm{d}x\ =\ +\infty .
\end{equation}
\end{thm}
Thus, we cannot expect that $u$ belongs to $W^{1,p}$ for $p>2$ for all time. In other words, the space $H^1=W^{1,2}$ is optimal for the equation \eqref{RB_2}.

Before proving Theorem \ref{optimality}, let $u_0 \in H^s$ with $s$ big enough, and let 
$u$ be a solution of $\mathrm{rB}$ with $u(0,x) = u_0(x)$. The main quantity is the following integral
\begin{equation}\label{intg}
\int_\mathcal{K} u_x^2(T,x)\, g(u_x(T,x))\, \mathrm{d}x,
\end{equation}
where $T>0$ and $\mathcal{K}$ is a compact set. Using the change of variable $x=y(T,\xi)$, one gets
\begin{equation}
\int_\mathcal{K} u_x^2(T,x)\, g(u_x(T,x))\, \mathrm{d}x\ =\ \int_\mathcal{K'} q\, \sin^2 (v\, /\, 2)\, g(\tan (v\, /\, 2))\, \mathrm{d} \xi ,
\end{equation}
where $\mathcal{K'}$ is another compact set. From previous sections, the quantity $q$ is always bounded, 
which implies that if $g$ is bounded then \eqref{intg} is bounded. 
If $g$ is not bounded (see Theorem \ref{optimality}), then the quantity 
\eqref{intg} depends on the behaviour of the derivative $u_x$ at time $T$. 
The proof of Theorem \ref{optimality} is done by building $u(T,\cdot)$, 
such that the quantity \eqref{intg} is infinite. Then, we use a backward system to go back 
in time and find a Lipschitz initial datum. 
\subsubsection*{Proof of Theorem \ref{optimality}:} 
Let $g(h) \eqdef \left[ \ln |h| \right]^\delta$ for $\delta>0$ and let $\bar{u}$ be a compactly supported odd function such that $\bar{u} \in \mathcal{C}^\infty(\mathds{R} / \{0 \})$ and for all $x \in ]0,\half [$ we have 
$$\bar{u}'(x)\ \eqdef\ - \frac{1}{\sqrt{x}}\, \left(- \ln  (x) \right)^{-\frac{1+\delta}{2}}.$$
It is clear that $\bar{u} \in H^1(\mathds{R})$ and 
\begin{equation}\label{ubar}
\int_{\mathcal{V}(0)} \bar{u}'(x)^2\, g(\bar{u}'(x))\, \mathrm{d}x\ =\ +\infty, \qquad \bar{u}'(x)\ \leqslant C,
\end{equation}
where $\mathcal{V}(0)$ denotes a neighbourhood of $0$.

The idea of the proof is to use a backward (in time) system such that $u(T,x)=\bar{u}(x)$. The initial datum $u_0$ is the unknown. To simplify the presentation, the conservative system \eqref{ODE} is used. With this system, we will obtain a local (in time) Oleinik inequality, which is enough for our construction. A similar proof can be used with the dissipative system \eqref{ODE2} with a global Oleinik inequality. The built solution in the interval $[0,T[$ is Lipchitz, so both systems \eqref{ODE}, \eqref{ODE2} yield the same solution.

In order to build $u_0$, we use the forward existence proof given in Section \ref{Global weak solutions}. One can use the change of variable $t \to -t$. The conservative system \eqref{ODE} becomes then
\begin{subequations}\label{systback}
\begin{align}
y_t\ &=\ -\, u,  &y(-T,\xi)\ =\ \bar{y}(\xi), \\ 
u_t\ &=\  \ell^2\, P_x,  &u(-T,\xi)\ =\ \bar{u}(\bar{y}(\xi)), \\ 
v_t\ &=\ P\left(1+\cos(v)\right)\ +\ \sin^2({v}/{2}),  &v(-T,\xi)\ =\ 2\/\arctan\!\left(
\bar{u}' \left(\bar{y}(\xi)\right)\right), \label{systbackc} \\ 
q_t\ &=\ -\, q\left(\half-P\right) \sin(v),   & q(-T,\xi)\ =\ 1,
\end{align}
\end{subequations}
where $t \in [-T,0]$ and $\bar{y}$ is defined as in \eqref{y}, replacing $u_0$ by $\bar{u}$.

The proof of a local existence of solutions can be done as in Section \ref{Global weak solutions}. Due to the change of variable $t \to -t$, the Oleinik inequality becomes
\begin{equation}
u_x(t,x)\ \geqslant\ -2/(t+T) 
\end{equation}
for $t>-T$ and $t$ close enough to $-T$. The proof of this Oleinik inequality proceeds as in Section \ref{Global weak solutions} using the equation \eqref{systbackc}, which implies that the derivative of the solution is bounded from below.
As in Remark \ref{Riccati_Oleinik}, since $\bar{u}'=u_x(-T,\cdot) \leqslant C$, the derivative of the solution remains bounded from above for $t>-T$ and $t$ close enough to $-T$.
Taking $T>0$ small so the solution is Lipschitz until $t=0$, and thus 
$$\int_\mathds{R} u_x(0,x)^2\, g(u_x(0,x)')\, \mathrm{d}x\ <\ +\infty .$$  
The result follows directly by using the change of variable $t \to -t$. \qed
\begin{rem} 
\begin{enumerate}
\item The optimality given in Theorem \ref{optimality} is also true for the Camassa--Holm equation.
\item \citet{XinZhang2000} have proved that the Camassa--Holm equation admits dissipative solutions that satisfy 
\begin{equation}\label{XZinq}
\int_0^T \int_{|x|\, \leqslant\, R} |u_x(t,x)|^p\, \mathrm{d}x\, \mathrm{d}t\ < +\infty \quad \forall T>0, R>0, p<3.
\end{equation}
This result can also be proven for the rB equation.
\item Theorem \ref{optimality} does not contradict with \eqref{XZinq}.
Theorem \ref{optimality} shows that the function 
$$t \mapsto \int_{|x|\, \leqslant\, R} |u_x(t,x)|^p\, \mathrm{d}x$$ 
does not necessarily belong to $L^\infty_{loc}([0,+\infty))$. However, the inequality \eqref{XZinq} shows that this function belongs to $L^1_{loc}([0,+\infty))$ if $p<3$.
\end{enumerate}
\end{rem}

\section{Conclusion and discussion}
In this paper, we have studied a regularisation of the inviscid Burgers equation \eqref{RB_2}. 
For a smooth initial datum, the regularised equation \eqref{RB_2} has a unique smooth solution 
locally in time. After the blow-up time, the solution is no longer unique, nor smooth. At least 
two types of solutions exist: conservative and dissipative solutions. We find that the {built}
dissipative solutions are more interesting because they satisfy an Oleinik inequality \eqref{Oleinik}, 
which plays an important role in showing that solutions converge {(up to a subsequence)} when 
$\ell \to 0$ and when $\ell \to \infty$ ($\ell$ the regularising positive parameter). Before the 
appearance of singularities, 
the limit when $\ell \to 0$ (respectively $\ell \to \infty$) is a smooth solution of the inviscid 
Burgers (resp. the Hunter--Saxton) equation. 
After {the breakdown time}, it remains open to determine whether the Burgers (resp. the 
Hunter--Saxton) equation holds in the limit without a remaining forcing term.



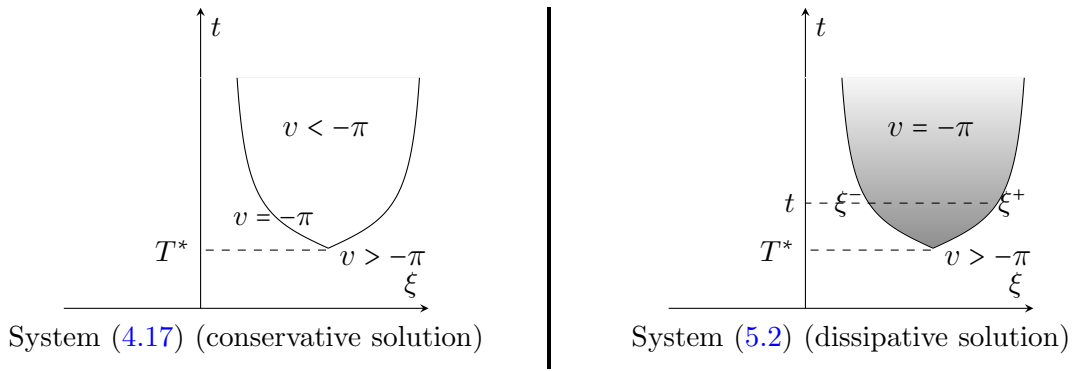
\begin{figure}[hbt!]
\begin{multicols}{2}
\begin{tikzpicture}[scale=0.7] 
\begin{axis}[
   axis lines=middle,
   xlabel={$\xi$},
   ylabel={$t$},
   xtick={\empty},
   ytick={\empty},
    domain=0.3:5,
   ymin=0,
   ymax=10,
   xmin=-3,
   xmax=5,
   clip=false,
]
\draw (2.75,6) node[] {\large $v < -\pi$};
\draw (1.6,3) node[] {$v = -\pi$};
\draw (4,1.7) node[] {\large $v >-\pi$};
\node at (1,-1) {System \eqref{ODE} (conservative solution)};

\draw[dashed]  (11/4,31/16) -- (0,31/16) node[anchor= east] {$T^*$}; ;
\addplot [mark max,black,name path=B, domain=0.8:4.8, samples=300] plot {bicuadratic(x)};

\Epigraph{0}{0}{0}{0}{B}

\end{axis}

\end{tikzpicture}

\begin{tikzpicture}[scale=0.7]
\begin{axis}[
   axis lines=middle,
   xlabel={$\xi$},
   ylabel={$t$},
   xtick={\empty},
   ytick={\empty},
    domain=0.3:5,
   ymin=0,
   ymax=10,
   xmin=-3,
   xmax=5,
   clip=false,
]
\draw[dashed]  (11/4,31/16) -- (0,31/16) node[anchor= east] {$T^*$}; ;
\draw[dashed]  (4,3.5) -- (0,3.5) node[anchor= east] {$t$};
\draw (4,3.5) node[right] {$\xi^+$};
\draw (1.5,3.5) node[left] {$\xi^-$};

\draw (2.75,6) node[] {\large $v=-\pi$};
\draw (4,1.7) node[] {\large $v >-\pi$};
\node at (1,-1) {System \eqref{ODE3} (dissipative solution)};

\addplot [mark max,black,name path=B, domain=0.8:4.8, samples=300] plot {bicuadratic(x)};
\DrawEpigraph{0}{0}{0}{0}{B}

\end{axis}
\end{tikzpicture}

\end{multicols}
\caption{Regions where $v=-\pi$.}
\label{fig:1}
\end{figure}

As shown above, the major difference between the conservative system 
\eqref{ODE} and the dissipative system \eqref{ODE3} is that the system 
\eqref{ODE} allows $v$ to cross the value $-\pi$, causing a jump of 
$u_x$ from $-\infty$ to $+\infty$ (see eq. \eqref{u_x}), which implies \eqref{AntiOl}, thence the loss of the Oleinik inequality (Remark \ref{remAntiOl}).
But, the value $v=-\pi$ is a barrier that cannot be crossed for the system \eqref{ODE3}. 
It follows that if $v(t,\xi_0)=-\pi$ at a time $t$, then $v(\tau,\xi_0)=-\pi$ for all 
times $\tau\geqslant t$ (see figure \ref{fig:1}). This property is important to obtain the Oleinik inequality \eqref{Oleinik}, 
which yields the dissipation of the energy \eqref{lossenergy}.

The figure \ref{fig:1} shows the domains where $v =-\pi$  for the systems \eqref{ODE} 
and \eqref{ODE3}.\\

%

\section*{Acknowledgments}
This material is based upon work supported by the 
National Science Foundation under Grant Nos.\ DMS  1812609 and 2106534 (RLP).

%

\end{document}